\documentclass[11pt]{article}

\usepackage[a4paper,
            includeheadfoot,
            margin=2cm]{geometry} 
\usepackage[utf8]{inputenc}       
\usepackage{lmodern}              
\usepackage{amsmath}              
\usepackage{amssymb}              
\usepackage{amsthm}               
\usepackage{commath}              
\usepackage{xspace}               
\usepackage{booktabs}             
\usepackage{numname}              
\usepackage[overload]{empheq}     
\usepackage[inline]{enumitem}     
\usepackage[pdftex]{hyperref}     

\newlist{myenum*}{enumerate*}{1}
\setlist[myenum*]{label=\itshape\roman*\upshape)}

\newcommand{\ie}{\textit{i.e.}\xspace}
\renewcommand{\vec}[1]{\mathbf{#1}}
\newcommand{\imag}{\jmath}

\newcommand{\deriv}[2]{\frac{\partial#1}{\partial#2}}
\newcommand{\derivv}[2]{\frac{\partial^2#1}{\partial#2}}
\DeclareMathOperator{\Grad}{\vec{grad}}
\DeclareMathOperator{\Div}{div}
\DeclareMathOperator{\Trsm}{\mathcal{S}}
\DeclareMathOperator{\Ident}{\mathcal{I}}
\DeclareMathOperator{\OpA}{\mathcal{A}}
\DeclareMathOperator{\Op}{Op}
\DeclareMathOperator{\DtN}{DtN}
\DeclareMathOperator{\Eig}{Eig}

\title{Convergence of classical optimized non-overlapping Schwarz method
  for Helmholtz problems in closed domains}
\author{Nicolas Marsic and Herbert De Gersem}
\date{\scriptsize Technische Universit\"at Darmstadt,
  Institut f\"ur Teilchenbeschleunigung und Elektromagnetische Felder (TEMF)}

\begin{document}
\maketitle
\begin{abstract}
  In this paper we discuss the convergence of state-of-the-art
  optimized Schwarz transmission conditions for Helmholtz problems
  defined on closed domains
  (\ie setups which do not exhibit an outgoing wave condition),
  as commonly encountered when modeling cavities.
  In particular, the impact of back-propagating waves
  on the Dirichlet-to-Neumann map is analyzed.
  Afterwards, the performance of
  the well-established optimized 0\textsuperscript{th}-order,
  evanescent modes damping, optimized 2\textsuperscript{nd}-order and
  Pad\'e-localized square-root transmission conditions is discussed.
\end{abstract}

\section{Introduction}
It is well known that large-scale time-harmonic Helmholtz problems
are hard to solve because of
\begin{myenum*}
\item the pollution effect~\cite{Ihlenburg1995} and
\item the indefiniteness of the discretized operator~\cite{Ernst2012}.
\end{myenum*}
While the pollution effect can be alleviated
by using higher order discretization schemes~\cite{Ihlenburg1997},
the indefiniteness is an intrinsic property of time-harmonic wave problems,
at least with standard variational formulations~\cite{Moiola2014, Diwan2019},
and significantly limits the performance of classical iterative solvers,
such as the generalized minimal residual method (GMRES).
Of course, as an alternative to iterative algorithms,
direct solvers can be used.
However, because of the fill-in effect,
whose minimization is know to be a NP-complete problem~\cite{Yannakakis1981},
the amount of memory needed to treat large-scale systems
can become prohibitively high (see for instance~\cite{Marsic2018b}).

As an alternative to direct and (unpreconditioned) iterative methods
for solving large-scale, high-frequency time-harmonic Helmholtz problems,
domain decomposition (DD) algorithms, and
optimized Schwarz (OS) techniques~\cite{Despres1990, Boubendir2007,
  Gander2002, Boubendir2012} in particular,
have attracted a lot of attention during the last decades.
The key idea thereof is to:
\begin{myenum*}
\item decompose the computational domain into (possibly overlapping) subdomains,
  creating thus new subproblems;
\item solve each subproblem \emph{independently};
\item exchange data at the interfaces between the subdomains
  via an appropriate \emph{transmission operator} and
\item solve each subproblem again and iterate until convergence of the solution.
\end{myenum*}
Since all subproblems are solved independently,
domain decomposition methods are parallel by nature\footnote{It is also possible
  to solve the subproblems sequentially
  and to exchange data after each single solve.
  This family of DD methods are often referred to as sweeping algorithms,
  and offer some advantages, notably in terms of iteration count,
  which will not be further discussed in this work.
  More details can be found for instance in~\cite{Vion2014a}.}
and are thus very well suited for the treatment of large-scale problems.
Furthermore, as the subproblems are of reduced size, direct solvers can be used.
Let also note that DD methods are rarely used as a stand-alone solver,
but most of the time as a \emph{preconditioner} for a Krylov subspace method
such as GMRES.
The design of such preconditioners for time-harmonic Helmholtz problems
remains an active and challenging topic~\cite{Gander2019}.

The convergence rate of an OS scheme strongly depends on its
transmission operator.
It is well known that the optimal operator is the
Dirichlet-to-Neumann ($\DtN$) map of the problem~\cite{Dolean2015a}
(\ie the operator relating the trace of the unknown field
to its normal derivative at the interface between two subdomains).
However, the $\DtN$ map is rarely employed
as it is a \emph{non-local} operator which leads to
a numerically expensive scheme.
Therefore, in practice, \emph{local approximations} of the $\DtN$ map are used,
which lead to many different computational schemes~\cite{Despres1990,
  Boubendir2007, Gander2002, Boubendir2012}
(see section~\ref{sec:open} for more details).
To the best of our knowledge,
all OS techniques share a common drawback:
they ignore the impact of \emph{back-propagating waves}.
While this assumption is legitimate in many cases
(antenna arrays~\cite{Peng2011},
medical imaging reconstruction~\cite{Tournier2017}
or photonic waveguides~\cite{Marsic2016a} just to cite a few),
it becomes questionable when the geometry allows resonances
(even if the source does not oscillate exactly at a resonance frequency),
as found for instance in lasers~\cite{Saleh2007},
accelerator cavities~\cite{Ko2006}
or quantum electrodynamic devices~\cite{Marsic2018b}.

The objective of this work is to determine the effect of
back-propagating waves on the performance of four well-established
transmission operators:
the optimized 0\textsuperscript{th}-order operator~\cite{Despres1990} (OO0),
the evanescent modes damping operator~\cite{Boubendir2007} (EMDA),
the optimized 2\textsuperscript{nd}-order operator~\cite{Gander2002} (OO2) and
the Pad\'e-localized square-root operator~\cite{Boubendir2012} (PADE).
This paper is organized as follows.
In section~\ref{sec:problem} the Helmholtz problem
as well as the optimized Schwarz scheme are presented formally
on a simple cavity model problem exhibiting back-propagating waves.
The optimal transmission operator of this model problem is then determined
in section~\ref{sec:S}.
Afterwards, in section~\ref{sec:comparison},
the optimal transmission condition for an unbounded problem without
obstacle (\ie exhibiting no back-propagating waves)
is recalled and compared with the one computed in the previous section.
The well-established OO0, EMDA, OO2 and PADE operators are recalled
in section~\ref{sec:open},
and their performance is analyzed for a rectangular cavity problems.
Section~\ref{sec:numexp} shows some numerical experiments, which
validate the previous theoretical analysis and
extend it to circular (2D) and spherical (3D) cavities.
Finally, conclusions are drawn in section~\ref{sec:conclusion}.

\section{Model problem and Schwarz scheme}\label{sec:problem}
\begin{figure}[ht]
  \centering
  \includegraphics{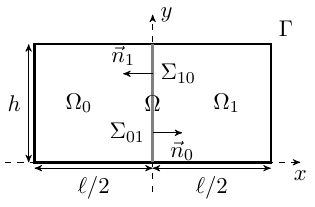}
  \caption{Domain $\Omega$ and its decomposition into $\Omega_1$ and $\Omega_2$.}
  \label{fig:domain}
\end{figure}

Let $\Omega$ be the two-dimensional domain $[-\ell/2, +\ell/2]\times[0, h]$
depicted in Figure~\ref{fig:domain}, and let $\Gamma$ be its boundary.
This domain is separated into two non-overlapping subdomains of equal size
$\Omega_0=[-\ell/2,0]\times[0, h]$ and
$\Omega_1=[0,+\ell/2]\times[0, h]$.
This splitting has introduced a new artificial boundary on each subdomain:
we denote by $\Sigma_{01}$ the artificial boundary of $\Omega_0$ and
by $\Sigma_{10}$ the artificial boundary of $\Omega_1$.
Furthermore,
$\vec{n}_i$ denotes the outwardly oriented unit vector normal to $\Sigma_{ij}$.

Let us solve the following Helmholtz problem on $\Omega$:
\begin{subequations}
  \label{eq:helmholtz}
  \begin{align}[left = \empheqlbrace]
    \Div{\Grad{p}} + k^2p & = g & \text{on}~\Omega,\label{eq:helmholtz:omega}\\
    p                     & = 0 & \text{on}~\Gamma,\label{eq:helmholtz:gamma}
  \end{align}
\end{subequations}
where $p(x, y)$ is the unknown function, $g(x, y)$ is a known source term
and $k\in\mathbb{R}$ is the fixed wavenumber of the Helmholtz problem.
Because of its boundary condition, it is obvious that~\eqref{eq:helmholtz}
models a \emph{cavity problem} strongly exhibiting forward-
\emph{and back-propagating waves}.
It is important to stress that for this problem to be well-defined,
we must assume that $k^2$ is not an eigenvalue of~\eqref{eq:helmholtz}.

Let us now set up the following optimized non-overlapping Schwarz method,
indexed by $n$, to solve the cavity Helmholtz problem~\eqref{eq:helmholtz}:
\begin{align*}[left = \empheqlbrace]
  \Div{\Grad{p_0^{n+1}}} + k^2p_0^{n+1}
  & = g
  & \text{on}~\Omega_0, \\
  \vec{n}_0\cdot\Grad{p_0^{n+1}} + \Trsm{(p_0^{n+1})}
  & =  \vec{n}_0\cdot\Grad{p_1^{n}} + \Trsm{(p_1^{n})}
  & \text{on}~\Sigma_{01}, \\
  p_0^{n+1}
  & = 0
  & \text{on}~\Gamma, \\
  \nonumber\\
  \Div{\Grad{p_1^{n+1}}} + k^2p_1^{n+1}
  & = g
  & \text{on}~\Omega_1, \\
  \vec{n}_1\cdot\Grad{p_1^{n+1}} + \Trsm{(p_1^{n+1})}
  & =  \vec{n}_1\cdot\Grad{p_0^{n}} + \Trsm{(p_0^{n})}
  & \text{on}~\Sigma_{10}, \\
  p_1^{n+1}
  & = 0
  & \text{on}~\Gamma,
\end{align*}
where $\Trsm$ is the transmission operator of the optimized Schwarz algorithm
and $p_i^n(x, y)$ is the solution of the iterative procedure
at iteration $n$ and on domain $\Omega_i$.
Once the Schwarz algorithm has converged,
the solution $p(x, y)$ of the original problem~\eqref{eq:helmholtz}
is recovered by concatenating the solutions $p_0(x, y)$ and $p_1(x, y)$.
Since the domains do not overlap, we furthermore have that
$\vec{n}_0 = -\vec{n}_1$.
Therefore, the system of equations becomes:
\begin{subequations}
  \label{eq:ddm}
  \begin{align}[left = \empheqlbrace]
    \Div{\Grad{p_0^{n+1}}} + k^2p_0^{n+1}
    & = g
    & \text{on}~\Omega_0, \\
    +\vec{n}_0\cdot\Grad{p_0^{n+1}} + \Trsm{(p_0^{n+1})}
    & = +\vec{n}_0\cdot\Grad{p_1^{n}} + \Trsm{(p_1^{n})}
    & \text{on}~\Sigma_{01}, \\
    p_0^{n+1}
    & = 0
    & \text{on}~\Gamma,  \label{eq:ddm:bc0} \\
    \nonumber\\
    \Div{\Grad{p_1^{n+1}}} + k^2p_1^{n+1}
    & = g
    & \text{on}~\Omega_1, \\
    -\vec{n}_0\cdot\Grad{p_1^{n+1}} + \Trsm{(p_1^{n+1})}
    & = -\vec{n}_0\cdot\Grad{p_0^{n}} + \Trsm{(p_0^{n})}
    & \text{on}~\Sigma_{10}, \\
    p_1^{n+1}
    & = 0
    & \text{on}~\Gamma. \label{eq:ddm:bc1}
  \end{align}
\end{subequations}

\section{Optimal transmission operator for the cavity problem}\label{sec:S}
In order to further simplify the problem,
let us now assume that the source term $g$ is zero.
Obviously, by not imposing a source in our problem
the solution $p(x,y)$ is trivially $p=0$ since $k$ is not an eigenvalue.
This however does not jeopardize the generality of
the results derived in this section.

Let us start by taking
the sine Fourier series of $p_i^n(x,y)$ along the $y$-axis:
\begin{equation}
  \label{eq:fourier}
  p_i^n(x,y) = \sum_{s\in\mathbb{S}}\widehat{p}_i^n(x, s)\sin(sy),
\end{equation}
where the functions $\widehat{p}_i^n(x, s)$ are the Fourier coefficients
and where $s$ is the Fourier variable, whose values are restricted to the set
\begin{equation}
  \label{eq:fourier:set}
  \mathbb{S} =
  \Big\{s\in\mathbb{R}~\Big|~s=m\frac{\pi}{h}, \forall m\in\mathbb{N}_0\Big\}.
\end{equation}
Indeed, by restricting $s$ to the set $\mathbb{S}$, the boundary conditions
\begin{align*}[left = \empheqlbrace]
  p^n_i(x, 0) & = 0 & \forall x\in\left[-\frac{\ell}{2}, +\frac{\ell}{2}\right],
  \\
  p^n_i(x, h) & = 0 & \forall x\in\left[-\frac{\ell}{2}, +\frac{\ell}{2}\right],
\end{align*}
are automatically satisfied.
Then, by exploiting decomposition~\eqref{eq:fourier},
the partial differential equation~\eqref{eq:ddm}
becomes the following ordinary differential equation (ODE):
\begin{subequations}
  \label{eq:ddm:fourier}
  \begin{align}[left = \empheqlbrace]
    \label{eq:ddm:fourier:ode:0}
    \derivv{\widehat{p}_0^{n+1}}{x^2} + (k^2-s^2)\widehat{p}_0^{n+1}
    & = 0
    & \forall x\in\left[-\frac{\ell}{2}, 0\right]~\text{and}~\forall s\in\mathbb{S}, \\
    \label{eq:ddm:fourier:trsm:0}
    +\deriv{\widehat{p}_0^{n+1}}{x}+\lambda\,\widehat{p}_0^{n+1}
    & = +\deriv{\widehat{p}_1^{n}}{x}+\lambda\,\widehat{p}_1^{n}
    & \text{on}~x=0~\text{and}~\forall s\in\mathbb{S}, \\
    \label{eq:ddm:fourier:bnd:0}
    \widehat{p}_0^{n+1}
    & = 0
    & \text{on}~x=-\frac{\ell}{2}~\text{and}~\forall s\in\mathbb{S}, \\
    \nonumber \\
    \label{eq:ddm:fourier:ode:1}
    \derivv{\widehat{p}_1^{n+1}}{x^2} + (k^2-s^2)\widehat{p}_1^{n+1}
    & = 0
    & \forall x\in\left[0, +\frac{\ell}{2}\right]~\text{and}~\forall s\in\mathbb{S}, \\
    \label{eq:ddm:fourier:trsm:1}
    -\deriv{\widehat{p}_1^{n+1}}{x}+\lambda\,\widehat{p}_1^{n+1}
    & = -\deriv{\widehat{p}_0^{n}}{x}+\lambda\,\widehat{p}_0^{n}
    & \text{on}~x=0~\text{and}~\forall s\in\mathbb{S}, \\
    \label{eq:ddm:fourier:bnd:1}
    \widehat{p}_1^{n+1}
    & = 0
    & \text{on}~x=+\frac{\ell}{2}~\text{and}~\forall s\in\mathbb{S},
  \end{align}
\end{subequations}
where $\lambda$ is the symbol of $\Trsm$.
Furthermore, and for simplicity, let us define $P_i^n(s)$ as:
\begin{equation}
  \label{eq:ddm:fourier:bnd:x0}
  P_i^n(s) = \widehat{p}_i^n(0, s).
\end{equation}

In order to find the best symbol $\lambda$,
we need to determine the convergence radius
of the iterative scheme~\eqref{eq:ddm:fourier}.
This objective can be achieved by:
\begin{enumerate}
\item deriving the fundamental solutions of~\eqref{eq:ddm:fourier:ode:0}
  and~\eqref{eq:ddm:fourier:ode:1};
  \label{enum:methodology:fundamental}

\item fixing the integration constants
  with the boundary conditions~\eqref{eq:ddm:fourier:bnd:0} and
  \eqref{eq:ddm:fourier:bnd:1} and the definition~\eqref{eq:ddm:fourier:bnd:x0};
  \label{enum:methodology:cst}

\item determining the solutions of~\eqref{eq:ddm:fourier:ode:0}
  and~\eqref{eq:ddm:fourier:ode:1}
  from the expressions found in steps~\ref{enum:methodology:fundamental}
  and~\ref{enum:methodology:cst};
  \label{enum:methodology:particular}

\item computing $\displaystyle\deriv{\widehat{p}_i^n(x, s)}{x}$ at $x=0$
  from the solutions $\widehat{p}_i^n(x,s)$ found
  in step~\ref{enum:methodology:particular} and
  \label{enum:methodology:derivative}

\item simplifying the transmission conditions~\eqref{eq:ddm:fourier:trsm:0}
  and~\eqref{eq:ddm:fourier:trsm:1}
  with the expressions found in steps~\ref{enum:methodology:particular}
  and~\ref{enum:methodology:derivative}.
\end{enumerate}
Let us note that
this methodology is the same as the one followed in~\cite{Gander2002}
for the Helmholtz problem in unbounded domains.

\subsection{Fundamental solutions for the case $s^2\neq{}k^2$}
The ODEs~\eqref{eq:ddm:fourier:ode:0} and~\eqref{eq:ddm:fourier:ode:1}
are nothing but a one-dimensional Helmholtz problem with wavenumber $k^2-s^2$.
Therefore, \emph{by assuming $s^2\neq{}k^2$}, the fundamental solutions are:
\begin{subequations}
  \label{eq:ddm:fourier:ode:fundamental}
  \begin{align}[left = \empheqlbrace]
    \label{eq:ddm:fourier:ode:fundamental:0}
    \widehat{p}_0^{n+1}(x,s)
    & = A_0~\exp\mathopen\Big[+\alpha(s)x\mathclose\Big] +
        B_0~\exp\mathopen\Big[-\alpha(s)x\mathclose\Big]
    & \forall x\leq0, \forall s\in\mathbb{S}, s^2\neq{}k^2, \\
    \label{eq:ddm:fourier:ode:fundamental:1}
    \widehat{p}_1^{n+1}(x,s)
    & = A_1~\exp\mathopen\Big[+\alpha(s)x\mathclose\Big] +
        B_1~\exp\mathopen\Big[-\alpha(s)x\mathclose\Big]
    & \forall x\geq0, \forall s\in\mathbb{S}, s^2\neq{}k^2,
  \end{align}
\end{subequations}
where $A_0$, $A_1$, $B_0$ and $B_1$ are integration constants,
and where
\begin{subequations}
  \label{eq:ddm:alpha}
  \begin{align}[left = {\alpha(s) = \empheqlbrace}]
    -\imag & \sqrt{k^2-s^2} & \text{if}~s^2\leq k^2, \\
           & \sqrt{s^2-k^2} & \text{if}~s^2\geq k^2,
  \end{align}
\end{subequations}
with $\imag$ the imaginary unit.
In what follows, only the case $s^2\neq{}k^2$ is discussed.
The alternative $s^2=k^2$ is addressed in section~\ref{sec:S:s=k}.

\subsection{Integration constants for the case $s^2\neq{}k^2$}\label{sec:S:cst}
Let us start by imposing the boundary conditions~\eqref{eq:ddm:fourier:bnd:x0}
and~\eqref{eq:ddm:fourier:bnd:0}.
By inserting them into~\eqref{eq:ddm:fourier:ode:fundamental:0},
we have for all $s\in\mathbb{S}$ and $s^2\neq{}k^2$:
\begin{align*}[left = \empheqlbrace]
  \widehat{p}_0^{n+1}(0, s)
  & =  P_0^{n+1}(s)
  &
  & \Longleftrightarrow
  & A_0 + B_0
  & =  P_0^{n+1}(s),
  \\
  \widehat{p}_0^{n+1}\left(-\frac{\ell}{2}, s\right)
  & = 0
  &
  & \Longleftrightarrow
  & A_0~\exp\mathopen{}\left[+\alpha(s)\,\frac{-\ell}{2}\mathclose{}\right] +
    B_0~\exp\mathopen{}\left[-\alpha(s)\,\frac{-\ell}{2}\mathclose{}\right]
  & = 0.
\end{align*}
Thus, it follows that:
\begin{equation*}
  \begin{array}{cl}
    & \left\{
      \begin{array}{r@{\,}l}
        A_0
        & = P_0^{n+1}(s) - B_0, \\
        B_0
        & = \displaystyle\frac%
          {-P_0^{n+1}(s)\,%
           \exp\mathopen{}\left[+\alpha(s)\,\frac{-\ell}{2}\mathclose{}\right]}%
          {\exp\mathopen{}\left[-\alpha(s)\,\frac{-\ell}{2}\mathclose{}\right]-%
           \exp\mathopen{}\left[+\alpha(s)\,\frac{-\ell}{2}\mathclose{}\right]},
      \end{array}
      \right. \\
    \\
    \Longleftrightarrow
    & \left\{
      \begin{array}{r@{\,}l}
        A_0
        & = +P_0^{n+1}(s)
          \left[
          1
          +
          \displaystyle\frac%
          {\exp\mathopen{}\left[-\alpha(s)\,\frac{\ell}{2}\mathclose{}\right]}%
          {\exp\mathopen{}\left[+\alpha(s)\,\frac{\ell}{2}\mathclose{}\right]-%
           \exp\mathopen{}\left[-\alpha(s)\,\frac{\ell}{2}\mathclose{}\right]}
          \right], \\
        \\
        B_0
        & = -P_0^{n+1}(s)\,%
          \displaystyle\frac%
          {\exp\mathopen{}\left[-\alpha(s)\,\frac{\ell}{2}\mathclose{}\right]}%
          {\exp\mathopen{}\left[+\alpha(s)\,\frac{\ell}{2}\mathclose{}\right]-%
           \exp\mathopen{}\left[-\alpha(s)\,\frac{\ell}{2}\mathclose{}\right]},
      \end{array}
      \right.
  \end{array}
\end{equation*}
and
\begin{subequations}
  \label{eq:ddm:fourier:ode:BC:0}
  \begin{align}[left = \empheqlbrace]
    A_0
    & = +P_0^{n+1}(s)\,%
      \frac{\exp\mathopen{}\left[+\alpha(s)\,\frac{\ell}{2}\mathclose{}\right]}%
           {\exp\mathopen{}\left[+\alpha(s)\,\frac{\ell}{2}\mathclose{}\right]-%
            \exp\mathopen{}\left[-\alpha(s)\,\frac{\ell}{2}\mathclose{}\right]},
    \\
    B_0
    & = -P_0^{n+1}(s)\,%
      \frac{\exp\mathopen{}\left[-\alpha(s)\,\frac{\ell}{2}\mathclose{}\right]}%
           {\exp\mathopen{}\left[+\alpha(s)\,\frac{\ell}{2}\mathclose{}\right]-%
            \exp\mathopen{}\left[-\alpha(s)\,\frac{\ell}{2}\mathclose{}\right]}.
  \end{align}
\end{subequations}

Similarly, the integration constants
of equation~\eqref{eq:ddm:fourier:ode:fundamental:1}
are found by inserting~\eqref{eq:ddm:fourier:bnd:x0}
and~\eqref{eq:ddm:fourier:bnd:1} into~\eqref{eq:ddm:fourier:ode:fundamental:1}:
\begin{align*}[left = \empheqlbrace]
  \widehat{p}_1^{n+1}(0, s)
  & =  P_1^{n+1}(s)
  &
  & \Longleftrightarrow
  & A_1 + B_1
  & =  P_1^{n+1}(s),
  \\
  \widehat{p}_1^{n+1}\left(+\frac{\ell}{2}, s\right)
  & = 0
  &
  & \Longleftrightarrow
  & A_1~\exp\mathopen{}\left[+\alpha(s)\,\frac{+\ell}{2}\mathclose{}\right] +
    B_1~\exp\mathopen{}\left[-\alpha(s)\,\frac{+\ell}{2}\mathclose{}\right]
  & = 0,
\end{align*}
which leads to
\begin{subequations}
  \label{eq:ddm:fourier:ode:BC:1}
  \begin{align}[left = \empheqlbrace]
    A_1
    & = +P_1^{n+1}(s)\,%
      \frac{\exp\mathopen{}\left[-\alpha(s)\,\frac{\ell}{2}\mathclose{}\right]}%
           {\exp\mathopen{}\left[-\alpha(s)\,\frac{\ell}{2}\mathclose{}\right]-%
            \exp\mathopen{}\left[+\alpha(s)\,\frac{\ell}{2}\mathclose{}\right]},
    \\
    B_1
    & = -P_1^{n+1}(s)\,%
      \frac{\exp\mathopen{}\left[+\alpha(s)\,\frac{\ell}{2}\mathclose{}\right]}%
           {\exp\mathopen{}\left[-\alpha(s)\,\frac{\ell}{2}\mathclose{}\right]-%
            \exp\mathopen{}\left[+\alpha(s)\,\frac{\ell}{2}\mathclose{}\right]}.
  \end{align}
\end{subequations}

\subsection{Solutions for the case $s^2\neq{}k^2$}
The solutions of the ODEs~\eqref{eq:ddm:fourier:ode:0}
and~\eqref{eq:ddm:fourier:ode:1},
subjected to the boundary conditions~\eqref{eq:ddm:fourier:bnd:x0},
\eqref{eq:ddm:fourier:bnd:0} and~\eqref{eq:ddm:fourier:bnd:1},
are then obtained by combining~\eqref{eq:ddm:fourier:ode:fundamental},
\eqref{eq:ddm:fourier:ode:BC:0} and~\eqref{eq:ddm:fourier:ode:BC:1}:
\begin{align*}[left = \empheqlbrace]
  \widehat{p}_0^{n+1}(x,s)
  & =
    +P_0^{n+1}(s)\,%
    \frac{\exp\mathopen{}\left[+\alpha(s)\,\frac{\ell}{2}\mathclose{}\right]}%
         {\exp\mathopen{}\left[+\alpha(s)\,\frac{\ell}{2}\mathclose{}\right]-%
          \exp\mathopen{}\left[-\alpha(s)\,\frac{\ell}{2}\mathclose{}\right]}~%
    \exp\mathopen\Big[+\alpha(s)x\mathclose\Big] \\
  & \phantom{=}
    -P_0^{n+1}(s)\,%
    \frac{\exp\mathopen{}\left[-\alpha(s)\,\frac{\ell}{2}\mathclose{}\right]}%
         {\exp\mathopen{}\left[+\alpha(s)\,\frac{\ell}{2}\mathclose{}\right]-%
          \exp\mathopen{}\left[-\alpha(s)\,\frac{\ell}{2}\mathclose{}\right]}~%
    \exp\mathopen\Big[-\alpha(s)x\mathclose\Big]
  & \forall x\leq0, \forall s\in\mathbb{S}, s^2\neq{}k^2,
  \\
  \\
  \widehat{p}_1^{n+1}(x,s)
  & =
    +P_1^{n+1}(s)\,%
    \frac{\exp\mathopen{}\left[-\alpha(s)\,\frac{\ell}{2}\mathclose{}\right]}%
         {\exp\mathopen{}\left[-\alpha(s)\,\frac{\ell}{2}\mathclose{}\right]-%
          \exp\mathopen{}\left[+\alpha(s)\,\frac{\ell}{2}\mathclose{}\right]}~%
    \exp\mathopen\Big[+\alpha(s)x\mathclose\Big] \\
  & \phantom{=}
    -P_1^{n+1}(s)\,%
    \frac{\exp\mathopen{}\left[+\alpha(s)\,\frac{\ell}{2}\mathclose{}\right]}%
         {\exp\mathopen{}\left[-\alpha(s)\,\frac{\ell}{2}\mathclose{}\right]-%
          \exp\mathopen{}\left[+\alpha(s)\,\frac{\ell}{2}\mathclose{}\right]}~%
    \exp\mathopen\Big[-\alpha(s)x\mathclose\Big]
  & \forall x\geq0, \forall s\in\mathbb{S}, s^2\neq{}k^2.
\end{align*}
Furthermore,
by definition of the hyperbolic sine\footnote{We have that~\cite{Oldham2009}:
  $2\,\sinh(x) = \exp(+x)-\exp(-x)$.}, we have:
\begin{subequations}
  \label{eq:ddm:fourier:ode:particular}
  \begin{align}[left = \empheqlbrace]
    \widehat{p}_0^{n+1}(x,s)
    & = +P_0^{n+1}(s)\,%
      \frac%
      {\exp\mathopen{}\left[+\alpha(s)\,\frac{\ell}{2}\mathclose{}\right]}%
      {2\,\sinh\mathopen{}\left[+\alpha(s)\,\frac{\ell}{2}\mathclose{}\right]}
      \exp\mathopen\Big[+\alpha(s)x\mathclose\Big]
      \nonumber\\
    & \phantom{=}
      -P_0^{n+1}(s)\,%
      \frac%
      {\exp\mathopen{}\left[-\alpha(s)\,\frac{\ell}{2}\mathclose{}\right]}%
      {2\,\sinh\mathopen{}\left[+\alpha(s)\,\frac{\ell}{2}\mathclose{}\right]}
      \exp\mathopen\Big[-\alpha(s)x\mathclose\Big]
    & \forall x\leq0, \forall s\in\mathbb{S}, s^2\neq{}k^2,
    \\
    \nonumber\\
    \widehat{p}_1^{n+1}(x,s)
    & =
      +P_1^{n+1}(s)\,%
      \frac%
      {\exp\mathopen{}\left[-\alpha(s)\,\frac{\ell}{2}\mathclose{}\right]}%
      {2\,\sinh\mathopen{}\left[-\alpha(s)\,\frac{\ell}{2}\mathclose{}\right]}
      \exp\mathopen\Big[+\alpha(s)x\mathclose\Big]
      \nonumber\\
    & \phantom{=}
      -P_1^{n+1}(s)\,%
      \frac%
      {\exp\mathopen{}\left[+\alpha(s)\,\frac{\ell}{2}\mathclose{}\right]}%
      {2\,\sinh\mathopen{}\left[-\alpha(s)\,\frac{\ell}{2}\mathclose{}\right]}
      \exp\mathopen\Big[-\alpha(s)x\mathclose\Big]
    & \forall x\geq0, \forall s\in\mathbb{S}, s^2\neq{}k^2.
  \end{align}
\end{subequations}

\subsection{Normal derivatives for the case $s^2\neq{}k^2$}
Thanks to the solution of
equation~\eqref{eq:ddm:fourier:ode:particular},
it is now possible to compute the normal derivatives of
$\widehat{p}_i^{n+1}(x,s)$:
\begin{subequations}
  \label{eq:ddm:fourier:ode:derivative}
  \begin{align}[left = \empheqlbrace]
    \deriv{\widehat{p}_0^{n+1}}{x}(x,s)
    & = +\alpha(s)\,P_0^{n+1}(s)\,%
      \frac%
      {\exp\mathopen{}\left[+\alpha(s)\,\frac{\ell}{2}\mathclose{}\right]}%
      {2\,\sinh\mathopen{}\left[+\alpha(s)\,\frac{\ell}{2}\mathclose{}\right]}
      \exp\mathopen\Big[+\alpha(s)x\mathclose\Big]
      \nonumber\\
    & \phantom{=}
      +\alpha(s)\,P_0^{n+1}(s)\,%
      \frac%
      {\exp\mathopen{}\left[-\alpha(s)\,\frac{\ell}{2}\mathclose{}\right]}%
      {2\,\sinh\mathopen{}\left[+\alpha(s)\,\frac{\ell}{2}\mathclose{}\right]}
      \exp\mathopen\Big[-\alpha(s)x\mathclose\Big]
    & \forall x\leq0, \forall s\in\mathbb{S}, s^2\neq{}k^2,
    \\
    \nonumber\\
    \deriv{\widehat{p}_1^{n+1}}{x}(x,s)
    & =
      +\alpha(s)\,P_1^{n+1}(s)\,%
      \frac%
      {\exp\mathopen{}\left[-\alpha(s)\,\frac{\ell}{2}\mathclose{}\right]}%
      {2\,\sinh\mathopen{}\left[-\alpha(s)\,\frac{\ell}{2}\mathclose{}\right]}
      \exp\mathopen\Big[+\alpha(s)x\mathclose\Big]
      \nonumber\\
    & \phantom{=}
      +\alpha(s)\,P_1^{n+1}(s)\,%
      \frac%
      {\exp\mathopen{}\left[+\alpha(s)\,\frac{\ell}{2}\mathclose{}\right]}%
      {2\,\sinh\mathopen{}\left[-\alpha(s)\,\frac{\ell}{2}\mathclose{}\right]}
      \exp\mathopen\Big[-\alpha(s)x\mathclose\Big]
    & \forall x\geq0, \forall s\in\mathbb{S}, s^2\neq{}k^2.
  \end{align}
\end{subequations}
Moreover, by evaluating these derivatives at $x=0$, it follows that:
\begin{align*}[left = \empheqlbrace]
  \deriv{\widehat{p}_0^{n+1}}{x}(0,s)
  & = +\alpha(s)\,P_0^{n+1}(s)\,%
    \frac%
    {\exp\mathopen{}\left[+\alpha(s)\,\frac{\ell}{2}\mathclose{}\right]}%
    {2\,\sinh\mathopen{}\left[+\alpha(s)\,\frac{\ell}{2}\mathclose{}\right]}
    \nonumber\\
  & \phantom{=}
    +\alpha(s)\,P_0^{n+1}(s)\,%
    \frac%
    {\exp\mathopen{}\left[-\alpha(s)\,\frac{\ell}{2}\mathclose{}\right]}%
    {2\,\sinh\mathopen{}\left[+\alpha(s)\,\frac{\ell}{2}\mathclose{}\right]}
  & \forall s\in\mathbb{S}, s^2\neq{}k^2,
  \\
  \nonumber\\
  \deriv{\widehat{p}_1^{n+1}}{x}(0,s)
  & =
    +\alpha(s)\,P_1^{n+1}(s)\,%
    \frac%
    {\exp\mathopen{}\left[-\alpha(s)\,\frac{\ell}{2}\mathclose{}\right]}%
    {2\,\sinh\mathopen{}\left[-\alpha(s)\,\frac{\ell}{2}\mathclose{}\right]}
    \nonumber\\
  & \phantom{=}
    +\alpha(s)\,P_1^{n+1}(s)\,%
    \frac%
    {\exp\mathopen{}\left[+\alpha(s)\,\frac{\ell}{2}\mathclose{}\right]}%
    {2\,\sinh\mathopen{}\left[-\alpha(s)\,\frac{\ell}{2}\mathclose{}\right]}
  & \forall s\in\mathbb{S}, s^2\neq{}k^2,
\end{align*}
which can be further simplified into
\begin{subequations}
  \label{eq:ddm:fourier:ode:derivative:x0}
  \begin{align}[left = \empheqlbrace]
    \deriv{\widehat{p}_0^{n+1}}{x}(0,s)
    & = +\alpha(s)\,P_0^{n+1}(s)\,%
      \coth\mathopen{}\left[+\alpha(s)\,\frac{\ell}{2}\mathclose{}\right]
    & \forall s\in\mathbb{S}, s^2\neq{}k^2,
    \\
    \nonumber\\
    \deriv{\widehat{p}_1^{n+1}}{x}(0,s)
    & = +\alpha(s)\,P_1^{n+1}(s)\,%
      \coth\mathopen{}\left[-\alpha(s)\,\frac{\ell}{2}\mathclose{}\right]
    & \forall s\in\mathbb{S}, s^2\neq{}k^2,
  \end{align}
\end{subequations}
by exploiting the definitions of the hyperbolic cosine  and
hyperbolic cotangent\footnote{We have that~\cite{Oldham2009}:
  $2\,\cosh(x) = \exp(+x)+\exp(-x)$ and $\coth(x) = \cosh(x)/\sinh(x)$.}.

\subsection{Convergence radius for the case $s^2\neq{}k^2$}
With the normal derivative of $\widehat{p}_i^{n+1}(x,s)$ in hand,
it is now possible to simplify
the transmission conditions~\eqref{eq:ddm:fourier:trsm:0}
and~\eqref{eq:ddm:fourier:trsm:1}.
By combining them with~\eqref{eq:ddm:fourier:ode:derivative:x0}
and~\eqref{eq:ddm:fourier:bnd:x0},
and by exploiting the parity of $\coth(x)$\footnote{The hyperbolic cotangent
  is an odd function~\cite{Oldham2009}: $\coth(-x) = -\coth(x)$},
we have:
\begin{equation*}
  \begin{array}{cl}
    & \left\{
      \begin{array}{r@{\,}l@{\qquad}r}
        +\displaystyle\deriv{\widehat{p}_0^{n+1}}{x}(0,s)
        +\lambda\,\widehat{p}_0^{n+1}(0,s)
        & = +\displaystyle\deriv{\widehat{p}_1^{n}}{x}(0,s)
          +\lambda\,\widehat{p}_1^{n}(0,s)
        & \forall s\in\mathbb{S}, s^2\neq{}k^2,
        \\
        -\displaystyle\deriv{\widehat{p}_1^{n+1}}{x}(0,s)
        +\lambda\,\widehat{p}_1^{n+1}(0,s)
        & = -\displaystyle\deriv{\widehat{p}_0^{n}}{x}(0,s)
          +\lambda\,\widehat{p}_0^{n}(0,s)
        & \forall s\in\mathbb{S}, s^2\neq{}k^2,
      \end{array}
      \right. \\
    \\
    \Longleftrightarrow
    & \left\{
      \begin{array}{r@{\,}l@{\qquad}r}
        & \displaystyle
          +\alpha(s)\,P_0^{n+1}(s)\,%
          \coth\mathopen{}\left[+\alpha(s)\,\frac{\ell}{2}\mathclose{}\right]
          +\lambda\,P_0^{n+1}(s)
        \\
        =
        & \displaystyle
          +\alpha(s)\,P_1^n(s)\,%
          \coth\mathopen{}\left[-\alpha(s)\,\frac{\ell}{2}\mathclose{}\right]
          +\lambda\,P_1^{n}(s)
        & \forall s\in\mathbb{S}, s^2\neq{}k^2, \\
        \\
        & \displaystyle
          -\alpha(s)\,P_1^{n+1}(s)\,%
          \coth\mathopen{}\left[-\alpha(s)\,\frac{\ell}{2}\mathclose{}\right]
          +\lambda\,P_1^{n+1}(s)
        \\
        =
        & \displaystyle
          -\alpha(s)\,P_0^n(s)\,%
          \coth\mathopen{}\left[+\alpha(s)\,\frac{\ell}{2}\mathclose{}\right]
          +\lambda\,P_0^{n}(s)
        & \forall s\in\mathbb{S}, s^2\neq{}k^2,
      \end{array}
      \right. \\
    \\
    \Longleftrightarrow
    & \left\{
      \begin{array}{r@{\,}l@{\qquad}r}
        \displaystyle
        P_0^{n+1}(s)
        \Bigg\{
        \lambda
        +
        \alpha(s)
        \coth\mathopen{}\left[\alpha(s)\,\frac{\ell}{2}\mathclose{}\right]
        \Bigg\}
        & =
          \displaystyle
          P_1^n(s)
          \Bigg\{
          \lambda
          -
          \alpha(s)
          \coth\mathopen{}\left[\alpha(s)\,\frac{\ell}{2}\mathclose{}\right]
          \Bigg\}
        & \forall s\in\mathbb{S}, s^2\neq{}k^2, \\
        \\
        \displaystyle
        P_1^{n+1}(s)
        \Bigg\{
        \lambda
        +
        \alpha(s)\,
        \coth\mathopen{}\left[\alpha(s)\,\frac{\ell}{2}\mathclose{}\right]
        \Bigg\}
        & =
          \displaystyle
          P_0^n(s)
          \Bigg\{
          \lambda
          -
          \alpha(s)
          \coth\mathopen{}\left[\alpha(s)\,\frac{\ell}{2}\mathclose{}\right]
        \Bigg\}
        & \forall s\in\mathbb{S}, s^2\neq{}k^2.
      \end{array}
      \right.
  \end{array}
\end{equation*}
Furthermore, since the index $n$ is arbitrary,
we can further simplify the two last equations into:
\begin{align*}[left = \empheqlbrace]
  P_0^{n+1}(s) & = (\rho^\text{close}_\lambda)^2(s)P_0^{n-1}(s)
               & \forall s\in\mathbb{S}, s^2\neq{}k^2, \\
  P_1^{n+1}(s) & = (\rho^\text{close}_\lambda)^2(s)P_1^{n-1}(s)
               & \forall s\in\mathbb{S}, s^2\neq{}k^2,
\end{align*}
where the convergence radius $\rho^\text{close}_\lambda(s)$ is given by
\begin{equation}
  \label{eq:fourier:rho:sneqk}
  \rho^\text{close}_\lambda(s) = \frac
  {\lambda-\alpha(s)
    \coth\mathopen{}\left[
      \alpha(s)\,\frac{\ell}{2}\mathclose{}
    \right]}
  {\lambda+\alpha(s)
    \coth\mathopen{}\left[
      \alpha(s)\,\frac{\ell}{2}\mathclose{}
    \right]}.
\end{equation}
From this last equation, it is then obvious that
the convergence radius can be reduced to $\rho^\text{close}_\lambda(s)=0$
for all $s^2\neq{}k^2$ by selecting:
\begin{equation}
  \label{eq:ddm:lambda:optimal:sneqk}
  \lambda = \lambda^{\text{opt}}_{\text{close}}(s) =
  \alpha(s)\coth\mathopen{}\left[
    \alpha(s)\,\frac{\ell}{2}
  \mathclose{}\right].
\end{equation}

\subsection{Case $s^2=k^2$}\label{sec:S:s=k}
Let us now treat the situation where $s^2=k^2$.
In this case, the ODEs~\eqref{eq:ddm:fourier:ode:0}
and~\eqref{eq:ddm:fourier:ode:1} admit as fundamental solution
$\widehat{p}_i^{n+1}(x,0) = A_i^{n+1}x+B^{n+1}_i$.
Then, by following the same approach as in section~\ref{sec:S:cst},
it is found directly that
$A_0^{n+1} = P_0\,2/\ell$, $A_1^{n+1} = -P_0\,2/\ell$,
and $B_i^{n+1} = P_i$.
Therefore, the normal derivatives are obviously
$\deriv{\widehat{p}_i^{n+1}}{x}(x,0)=A_i^{n+1}$,
leading thus to a convergence radius of the OS scheme of the form:
\begin{equation}
  \label{eq:fourier:rho:s=k}
  \rho^\text{close}_\lambda(0) = \frac{\lambda-2/\ell}{\lambda+2/\ell}.
\end{equation}
This convergence radius can thus be reduced to $\rho^\text{close}_\lambda(0)=0$
by selecting
\begin{equation}
  \label{eq:ddm:lambda:optimal:s=k}
  \lambda = \lambda^{\text{opt}}_{\text{close}}(0) = 2/\ell.
\end{equation}

\subsection{Optimal operator}
By summarizing the results obtained in~\eqref{eq:fourier:rho:sneqk},
\eqref{eq:fourier:rho:s=k}, \eqref{eq:ddm:lambda:optimal:sneqk}
and \eqref{eq:ddm:lambda:optimal:s=k},
it follows that the optimal Schwarz operator
\begin{equation}
  \label{eq:close:S:opt}
  \Trsm^{\text{opt}}_{\text{close}} =
  \Op\mathopen{}\Bigg\{\lambda^{\text{opt}}_{\text{close}}\mathclose{}\Bigg\}
\end{equation}
has as symbol:
\begin{subequations}
  \label{eq:close:lambda:opt}
  \begin{align}[left = {\lambda^{\text{opt}}_{\text{close}}(s)=\empheqlbrace}]
    & \sqrt{k^2-s^2}
      \cot\mathopen{}\left[\sqrt{k^2-s^2}\frac{\ell}{2}\mathclose{}\right]
    & \text{if}~s^2<k^2,\\
    & 2/\ell
    & \text{if}~s^2=k^2,\\
    & \sqrt{s^2-k^2}
      \coth\mathopen{}\left[\sqrt{s^2-k^2}\frac{\ell}{2}\mathclose{}\right]
    & \text{if}~s^2>k^2,
  \end{align}
\end{subequations}
since $\coth\mathopen(\imag{}a\mathclose)=-\imag\cot\mathopen(a\mathclose)$
and $\cot\mathopen(-a\mathclose)=-\cot\mathopen(a\mathclose)$
$\forall a\in\mathbb{R}$~\cite{Oldham2009}.
Furthermore, the associated convergence radius is given by:
\begin{equation}
  \label{eq:fourier:rho}
  \rho^\text{close}_\lambda(s) =
  \frac{\lambda(s)-\lambda^{\text{opt}}_{\text{close}}(s)}
       {\lambda(s)+\lambda^{\text{opt}}_{\text{close}}(s)}.
\end{equation}
From this last equation, it is clear that $\rho^\text{close}_\lambda(s)=0$
if we select $\lambda = \lambda^{\text{opt}}_{\text{close}}$.

\section{Comparison between the optimal operators
  for cavity problems and unbounded problems without obstacles}
\label{sec:comparison}
Let us now consider the following \emph{unbounded} Helmholtz problem
without obstacles:
\begin{subequations}
  \label{eq:helmholtz:open}
  \begin{align}[left = \empheqlbrace]
    \Div{\Grad{p}} + k^2p & = g \quad\text{on}~\mathbb{R}^2,\\
    \lim_{r\to\infty}\sqrt{r}\left(\deriv{p}{r}-\imag{}kp\right) & = 0,
    \label{eq:helmholtz:open:rad}
  \end{align}
\end{subequations}
where $r^2 = x^2+y^2$.
In this case, it can be shown that the optimal transmission operator
$\Trsm^{\text{opt}}_{\text{open}}$
for solving this problem with an OS scheme writes~\cite{Boubendir2012}:
\begin{equation}
  \label{eq:open:S:opt}
  \Trsm^{\text{opt}}_{\text{open}} =
  \Op\mathopen{}\Big(\lambda^{\text{opt}}_{\text{open}}\mathclose{}\Big) =
  -\imag{}k\sqrt{1+\frac{\Div_\Sigma\Grad_\Sigma}{k^2}},
\end{equation}
where
\begin{equation}
  \label{eq:open:lambda:opt}
  \lambda^{\text{opt}}_{\text{open}} = -\imag{}k\sqrt{1-\frac{s^2}{k^2}}.
\end{equation}

\begin{figure}[ht]
  \centering
  \includegraphics{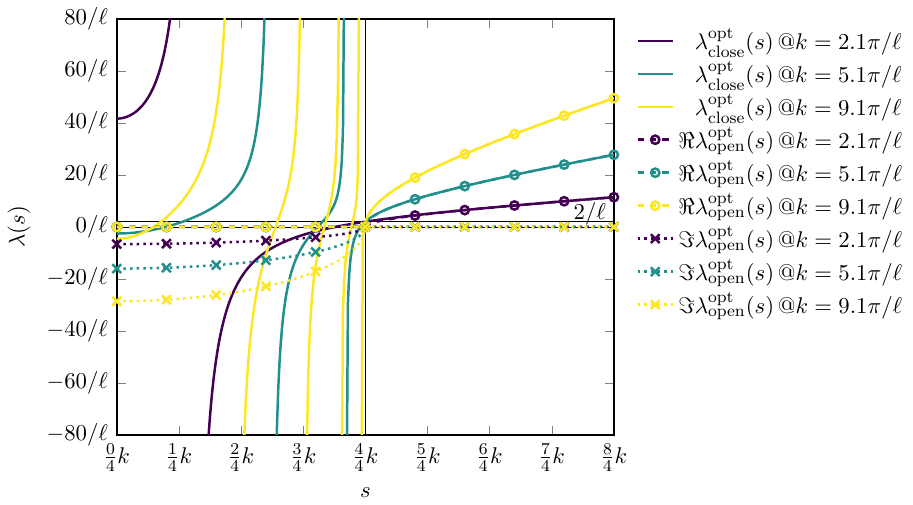}
  \caption{Graphs of
    $\lambda^{\text{opt}}_{\text{open}}(s)$ and
    $\lambda^{\text{opt}}_{\text{close}}(s)$
    for different values of $k$.}
  \label{fig:lambda:k}
\end{figure}

This non-local operator
is the keystone of the construction of the OO0, EMDA, OO2 and PADE
transmission operators.
In particular, the four aforementioned transmission operators
are nothing else but \emph{local approximations} of
$\Trsm^{\text{opt}}_{\text{open}}$.
Therefore, before studying the performance of OO0, EMDA, OO2 and PADE
for solving the \emph{cavity} problem~\eqref{eq:ddm},
let us first compare the two optimal operators
$\Trsm^{\text{opt}}_{\text{open}}$ and $\Trsm^{\text{opt}}_{\text{close}}$,
or, more precisely, their symbols.
By comparing~\eqref{eq:open:lambda:opt} and~\eqref{eq:close:lambda:opt},
it is easy to realize that
\begin{subequations}
  \label{eq:open:close}
  \begin{align}[left = {\lambda^{\text{opt}}_{\text{close}}(s)-\lambda^{\text{opt}}_{\text{open}}(s)=\empheqlbrace}]
    & \sqrt{k^2-s^2}
      \cot\mathopen{}\left[\sqrt{k^2-s^2}\frac{\ell}{2}\mathclose{}\right]
      +\imag{}\sqrt{k^2-s^2}
    & \text{if}~s^2<k^2, \label{eq:open:close:s<k}\\
    & 2/\ell
    & \text{if}~s^2=k^2,\\
    & \sqrt{s^2-k^2}
      \coth\mathopen{}\left[\sqrt{s^2-k^2}\frac{\ell}{2}\mathclose{}\right]
      -\sqrt{s^2-k^2}
    & \text{if}~s^2>k^2.
  \end{align}
\end{subequations}
Interestingly, by exploiting the definition of the hyperbolic
cotangent~\cite{Oldham2009}, the case $s^2>k^2$ can be further simplified into
\begin{align}
  \lambda^{\text{opt}}_{\text{close}}(s)-\lambda^{\text{opt}}_{\text{open}}(s)
  & = \sqrt{s^2-k^2}
    \coth\mathopen{}\left[\sqrt{s^2-k^2}\frac{\ell}{2}\mathclose{}\right]
    -\sqrt{s^2-k^2}\nonumber\\
  & = \sqrt{s^2-k^2}
      \left(\frac{\exp\mathopen{}\left(\ell\sqrt{s^2-k^2}\mathclose{}\right)+1}
                 {\exp\mathopen{}\left(\ell\sqrt{s^2-k^2}\mathclose{}\right)-1}
            - 1\right)\nonumber\\
  & = \frac{2}{\exp\mathopen{}\left(\ell\sqrt{s^2-k^2}\mathclose{}\right)-1}
      \sqrt{s^2-k^2}
  & \text{if}~s^2>k^2, \label{eq:open:close:s>k}
\end{align}
which yields to:
\begin{equation}
  \label{eq:open:close:lim}
  \lim_{s\to\infty}
  \lambda^{\text{opt}}_{\text{close}}(s)-\lambda^{\text{opt}}_{\text{open}}(s)
  = 0.
\end{equation}
In other words, for the case $s^2>k^2$,
the symbol $\lambda^{\text{opt}}_{\text{open}}(s)$
is converging towards $\lambda^{\text{opt}}_{\text{close}}(s)$ as $s$ grows.
Furthermore,
as the difference between those two symbols is decreasing exponentially,
\emph{$\lambda^{\text{opt}}_{\text{open}}(s)$
is an excellent approximation of $\lambda^{\text{opt}}_{\text{close}}(s)$
when $s^2>k^2$.}
For the case $s^2<k^2$, as the codomains of
$\lambda^{\text{opt}}_{\text{open}}(s)$ (which is purely imaginary)
and $\lambda^{\text{opt}}_{\text{close}}(s)$ (which is purely real)
do not match, the expression in~\eqref{eq:open:close:s<k}
cannot be further simplified.
Moreover, the poles and zeros of the two symbols are very different
for the case $s^2<k^2$, as summarized in Table~\ref{tab:open:close}.
For illustration purposes, the graphs of
$\lambda^{\text{opt}}_{\text{open}}$
and $\lambda^{\text{opt}}_{\text{close}}$
are depicted in Figure~\ref{fig:lambda:k} for different values of $k$
($\Re$ and $\Im$ are respectively denoting
the real and imaginary part functions).
\begin{table}[ht]
  \centering
  \begin{tabular}{ccc}
    \hline
    \hline
    Criterion
    & $\lambda^{\text{opt}}_{\text{open}}(s)$
    & $\lambda^{\text{opt}}_{\text{close}}(s)$ \\
    \cmidrule(lr){1-1}
    \cmidrule(lr){2-2}
    \cmidrule(lr){3-3}
    Codomain
    & $\imag\mathbb{R}^-$ (if $s^2\leq{}k^2$),
      $\mathbb{R}^+$ (if $s^2\geq{}k^2$)
    & $\mathbb{R}$ \\
    Zeros
    & one at $s=k$
    & many (if $s^2<k^2$), none (if $s^2>k^2$) \\
    Poles
    & none
    & many (if $s^2<k^2$), none (if $s^2>k^2$) \\
    Value at $s=0$
    & $0$
    & $2/\ell$ \\
    \hline
    \hline
  \end{tabular}
  \caption{Comparison between $\lambda^{\text{opt}}_{\text{open}}$ and
    $\lambda^{\text{opt}}_{\text{close}}$.}
  \label{tab:open:close}
\end{table}

The (dis)similarities between $\lambda^{\text{opt}}_{\text{open}}(s)$ and
$\lambda^{\text{opt}}_{\text{close}}(s)$
discussed above from a mathematical point of view
can also be given a more physical interpretation.
From the analysis carried out in section~\ref{sec:S},
we know that the Dirichlet-to-Neumann ($\DtN$) map
of problem~\eqref{eq:helmholtz} on $\Sigma_{10}$ (resp. $\Sigma_{01}$)
is given by $\Trsm^{\text{opt}}_{\text{close}}$.
Therefore, when approximating this $\DtN$ map
by $\Trsm^{\text{opt}}_{\text{open}}$,
we assume that beyond $\Sigma_{10}$ (resp. $\Sigma_{01}$)
the solution $p(x,y)$ will behave as a wave in open space.
This approximation makes sense for the component $\widehat{p}(x,s)$
such that $s^2>k^2$.
Indeed, in this case,
the solution $\widehat{p}(x,s)$ is an \emph{evanescent wave},
which decays as it goes away from $x=0$.
Therefore, it makes almost no difference whether the domain is closed or open
for these components
and the use of $\Trsm^{\text{opt}}_{\text{open}}$ is thus legitimate.
On the other hand, the components $\widehat{p}(x,s)$ such that $s^2<k^2$
\emph{corresponds to waves propagating}
in the subdomains $\Omega_0$ and $\Omega_1$.
Obviously, assuming that both $\Omega_0$ and $\Omega_1$ correspond
to unbounded domains is incorrect,
meaning that $\lambda^{\text{opt}}_{\text{open}}(s)$ cannot be used
as a good approximation of $\lambda^{\text{opt}}_{\text{close}}(s)$
in this latter case.

\section{Behavior of transmission operators optimized for unbounded problem used
  in a cavity configuration}
\label{sec:open}
Now that we have presented and compared the optimal operators
$\Trsm^{\text{opt}}_{\text{close}}$ and $\Trsm^{\text{opt}}_{\text{open}}$,
we can determine the performance of the
OO0, EMDA, OO2 and PADE operators
(which are nothing but approximations of  $\Trsm^{\text{opt}}_{\text{open}}$),
when used as an approximations of $\Trsm^{\text{opt}}_{\text{close}}$.
However, before starting this study,
let us recall shortly the four transmission operators.

\subsection{Recall of the OO0, EMDA, OO2 and PADE operators}

\subsubsection{Optimized zeroth-order optimized operator (OO0)}
The simplest approximate of $\Trsm^{\text{opt}}_{\text{open}}$
is simply a constant value~\cite{Despres1990},
which is selected from the constant term
of the Taylor expansion of $\lambda^{\text{opt}}_{\text{open}}$.
This leads to the so-called optimized 0\textsuperscript{th}-order operator (OO0),
whose symbol reads:
\begin{equation}
  \label{eq:open:lambda:oo0}
  \lambda_{\text{open}}^{\text{OO0}} = -\imag{}k.
\end{equation}

\subsubsection{Evanescent modes damping operator (EMDA)}
In order to further increase the performance of the OO0 operator,
a complexified wavenumber $k_\varepsilon$ can be introduced:
\begin{equation}
  \label{eq:open:keps}
  k_\varepsilon = (1 + \imag\varepsilon)k,
\end{equation}
where $\varepsilon$ is a positive real value.
This complexification leads then to the so-called
evanescent modes damping operator~\cite{Boubendir2007} (EMDA),
whose symbol reads:
\begin{equation}
  \label{eq:open:lambda:emda}
  \lambda_{\text{open}}^{\text{EMDA}} = -\imag{}k_\varepsilon.
\end{equation}

\subsubsection{Optimized second-order optimized operator (OO2)}
By pushing the Taylor approximation strategy further,
a second-order symbol can be designed:
\begin{equation}
  \label{eq:open:lambda:oo2}
  \lambda_{\text{open}}^{\text{OO2}} = a+bs^2,
\end{equation}
where $a$ and $b$ are two complex-valued constants,
chosen to optimize the convergence rate of the Schwarz scheme~\cite{Gander2002}.
Let us note that the optimal choice for $a$ and $b$
differs from the coefficients of the Taylor expansion.
The operator associated with this symbol is classically referred to as
the optimized second-order operator (OO2).

\subsubsection{Pad\'e-localized square-root transmission condition (PADE)}
As an alternative to the Taylor expansion, a Pad\'e decomposition
of the square-root symbol~\eqref{eq:open:lambda:opt} can be employed.
This strategy leads to the following approximation
with $N_p$ Pad\'e terms~\cite{Boubendir2012}:
\begin{equation}
  \label{eq:open:lambda:pade}
  \lambda_{\text{open}}^{\text{Pad\'e}}
  =
  -\imag kC_0
  -\imag k\sum_{p=1}^{N_p}
         \left(A_{p}\frac{s^2}{k_\varepsilon^2}\right)
         \left(1 + B_{p}\frac{s^2}{k_\varepsilon^2}\right)^{-1}.
\end{equation}
The coefficients $C_0$, $A_p$ and $B_p$ are given by
\begin{equation*}
C_0 = e^{\imag\xi/2}R_{N_p}\left(e^{-\imag \xi} - 1\right),
\quad
A_{p} = \frac{e^{-\imag\xi/2}a_{p}}
                {\big[1 + b_{p}(e^{-\imag \xi}-1)\big]^{2}},
\quad
B_{p} = \frac{e^{-\imag \xi}b_{p}}
                {1 + b_{p}(e^{-\imag \xi}-1)},
\end{equation*}
where:
\begin{itemize}
\item $\xi$ is a rotation angle of the branch cut of the square-root function
  and is usually taken as $\pi/4$;
\item $R_{N_p}(z)$ is the standard real-valued Pad\'e approximation of
  order $N_p$ of $\sqrt{1+z}$, that is
  \begin{equation*}
    R_{N_p}(z) = 1+\sum_{p=1}^{N_p} \frac{a_{p} z}{1 + b_{p} z};
  \end{equation*}
\item $a_p$ and $b_p$ are defined as
  \begin{equation*}
    a_{p} = \frac{2}{2N_p+1}\sin^{2}\mathopen{}\left(
      \frac{p\pi}{2N_p+1}
    \mathclose{}\right),
    \qquad
    b_{p} = \cos^{2}\mathopen{}\left(\frac{p\pi}{2N_p+1}\mathclose{}\right).
  \end{equation*}
\end{itemize}

\subsection{Transmission operators for the open problem
  as an approximation of the optimal transmission operator
  for the closed problem}
\label{sec:open:S}
Now that the OO0, EMDA, OO2 and PADE operators are recalled,
let us analyze their performance,
when the symbols $\lambda_{\text{open}}^{\text{OO0}}$,
$\lambda_{\text{open}}^{\text{EMDA}}$, $\lambda_{\text{open}}^{\text{OO2}}$
and $\lambda_{\text{open}}^{\text{Pad\'e}}$
are used to approximate $\lambda^{\text{opt}}_{\text{close}}$.
Because of the asymptotic behavior of
$\lambda^{\text{opt}}_{\text{close}}(s)-\lambda^{\text{opt}}_{\text{open}}(s)$
given in~\eqref{eq:open:close:lim},
we already know that OO0, EMDA, OO2 and PADE
will be good approximations of $\Trsm^{\text{opt}}_{\text{close}}$
(or, at least, as good as they were for $\Trsm^{\text{opt}}_{\text{open}}$)
for the \emph{evanescent modes}.

\subsubsection{Optimized zeroth-order operator (OO0)}
\begin{figure}[ht]
  \centering
  \includegraphics{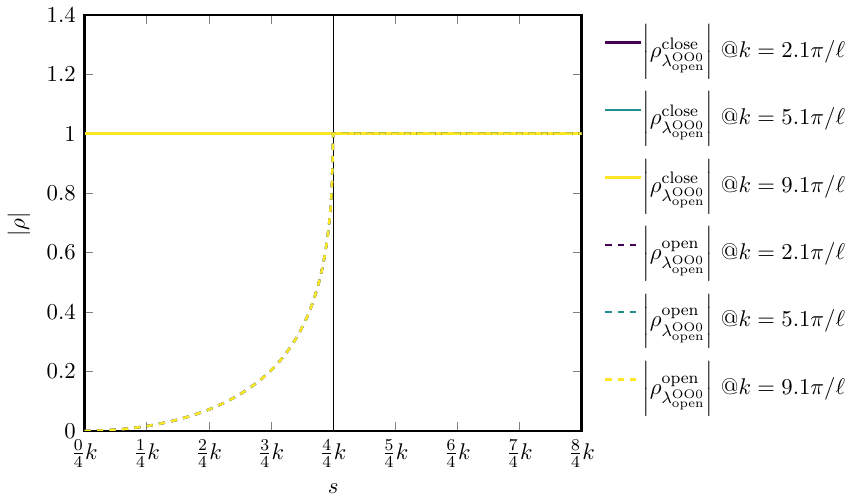}
  \caption{Graphs of
    $\abs{\rho^\text{close}_{\lambda_{\text{open}}^{\text{OO0}}}}$ and
    $\abs{\rho^\text{open}_{\lambda_{\text{open}}^{\text{OO0}}}}$
    for different values of $k$.}
  \label{fig:rho:oo0}
\end{figure}

In the case of the OO0 symbol~\eqref{eq:open:lambda:oo0},
and given the convergence radius $\rho^\text{close}_\lambda(s)$
of the Schwarz scheme~\eqref{eq:ddm},
we have that:
\begin{equation*}
  \rho^\text{close}_{\lambda_\text{open}^{\text{OO0}}}(s) =
  \frac{-\imag{}k-\lambda^{\text{opt}}_{\text{close}}(s)}
       {-\imag{}k+\lambda^{\text{opt}}_{\text{close}}(s)}.
\end{equation*}
Therefore, since $\lambda^{\text{opt}}_{\text{close}}(s)$ is real-valued,
the modulus of $\rho^\text{close}_{\lambda_\text{open}^{\text{OO0}}}(s)$
is then simply:
\begin{equation}
  \label{eq:close:rho:oo0}
  \abs{\rho^\text{close}_{\lambda_\text{open}^{\text{OO0}}}(s)}^2 = 1
  \qquad\forall{}s\in\mathbb{S}.
\end{equation}
This last result can be compared to the unbounded case,
where the modulus of the convergence radius
$\abs{\rho^\text{open}_{\lambda_\text{open}^{\text{OO0}}}(s)}$
is~\cite{Dolean2015a} (assuming no overlap):
\begin{subequations}
  \label{eq:open:rho:oo2}
  \begin{align}[left =
    {\abs{\rho^\text{open}_{\lambda_\text{open}^{\text{OO0}}}(s)}\empheqlbrace}]
    & < 1 & \text{if}~s^2<k^2,\\
    & = 1 & \text{if}~s^2\geq{}k^2.
  \end{align}
\end{subequations}
As expected, the behaviors of
$\abs{\rho^\text{close}_{\lambda_\text{open}^{\text{OO0}}}(s)}$
and
$\abs{\rho^\text{open}_{\lambda_\text{open}^{\text{OO0}}}(s)}$
are identical when $s^2>k^2$.
On the other hand, compared to the unbounded problem,
the performance of OO0 is significantly deteriorated
when solving a cavity problem.
For illustration purposes, the graphs of
$\abs{\rho^\text{close}_{\lambda_\text{open}^{\text{OO0}}}(s)}$
and
$\abs{\rho^\text{open}_{\lambda_\text{open}^{\text{OO0}}}(s)}$
are shown in Figure~\ref{fig:rho:oo0}.

\subsubsection{Evanescent modes damping operator (EMDA)}
In the case of the EMDA symbol~\eqref{eq:open:lambda:emda},
the convergence radius
$\rho^\text{close}_{\lambda_\text{open}^{\text{EMDA}}}(s)$
writes
\begin{equation*}
  \rho^\text{close}_{\lambda_\text{open}^{\text{EMDA}}}(s) =
  \frac{-\imag{}k-(\lambda^{\text{opt}}_{\text{close}}(s)-\varepsilon{}k)}
       {-\imag{}k+(\lambda^{\text{opt}}_{\text{close}}(s)+\varepsilon{}k)},
\end{equation*}
and its modulus is given by:
\begin{equation}
  \label{eq:close:rho:emda}
  \abs{\rho^\text{close}_{\lambda_\text{open}^{\text{EMDA}}}(s)}^2 =
  \frac{k^2+(\lambda^{\text{opt}}_{\text{close}}(s)-\varepsilon{}k)^2}
       {k^2+(\lambda^{\text{opt}}_{\text{close}}(s)+\varepsilon{}k)^2}.
\end{equation}
Depending on the values of $\lambda^{\text{opt}}_{\text{close}}(s)$, it is easy
to see that $\rho^\text{close}_{\lambda_\text{open}^{\text{EMDA}}}(s)$
exhibits the following properties:
\begin{enumerate}
\item $\abs{\rho^\text{close}_{\lambda_\text{open}^{\text{EMDA}}}(s)}^2\to1$
  if $\lambda^{\text{opt}}_{\text{close}}(s)\to\pm\infty$,
  $\lambda^{\text{opt}}_{\text{close}}(s)\to0$
  or $\varepsilon\to0$;
\item $\abs{\rho^\text{close}_{\lambda_\text{open}^{\text{EMDA}}}(s)}^2<1$
  if $\lambda^{\text{opt}}_{\text{close}}(s)>0$ and $\varepsilon>0$;
\item $\abs{\rho^\text{close}_{\lambda_\text{open}^{\text{EMDA}}}(s)}^2>1$
  if $\lambda^{\text{opt}}_{\text{close}}(s)<0$ and $\varepsilon>0$.
\end{enumerate}
Obviously, the two last results are inverted for $\varepsilon<0$.

\begin{figure}[ht]
  \centering
  \includegraphics{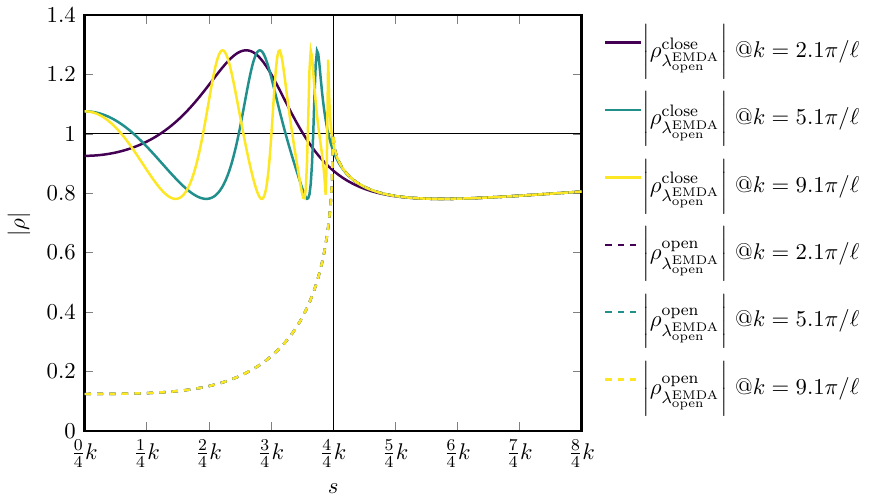}
  \caption{Graphs of
    $\abs{\rho^\text{close}_{\lambda_{\text{open}}^{\text{EMDA}}}}$ and
    $\abs{\rho^\text{open}_{\lambda_{\text{open}}^{\text{EMDA}}}}$
    for different values of $k$ and with $\varepsilon = 0.25$.}
  \label{fig:rho:emda}
\end{figure}

Let us now compare $\rho^\text{close}_{\lambda_{\text{open}}^{\text{EMDA}}}$
and
$\rho^\text{open}_{\lambda_{\text{open}}^{\text{EMDA}}}$
(\ie the convergence radius of the EMDA operator when used
in an OS scheme for solving the unbounded problem~\eqref{eq:helmholtz:open}).
As shown in Figure~\ref{fig:rho:emda} for a damping
coefficient of $\varepsilon=25\%$,
the difference between the two radii becomes
unnoticeable as $s$ grows (once $s^2>k^2$).
On the other hand,
we have that:
\begin{equation*}
  \max_s\abs{\rho^\text{open}_{\lambda_{\text{open}}^{\text{EMDA}}}(s)}<1~
  \text{(see~\cite{Boubendir2007})}
  \qquad
  \max_s\abs{\rho^\text{close}_{\lambda_{\text{open}}^{\text{EMDA}}}(s)}>1
  \qquad
  \forall{}s^2<k^2.
\end{equation*}
In other words,
the performance of EMDA is deteriorated
when passing from an unbounded wave problem without obstacle to a cavity one.

\subsubsection{Optimized second-order operator (OO2)}
Let us now focus on the OO2 symbol~\eqref{eq:open:lambda:oo2}.
In this case, the convergence radius
$\rho^\text{close}_{\lambda_\text{open}^{\text{OO2}}}(s)$
is given by
\begin{equation*}
  \rho^\text{close}_{\lambda_\text{open}^{\text{OO2}}}(s) =
  \frac{(a+bs^2)-\lambda^{\text{opt}}_{\text{close}}(s)}
       {(a+bs^2)+\lambda^{\text{opt}}_{\text{close}}(s)},
\end{equation*}
and its modulus reads:
\begin{equation}
  \label{eq:close:rho:oo2}
  \abs{\rho^\text{close}_{\lambda_\text{open}^{\text{OO2}}}(s)}^2 =
  \frac{\abs{a+bs^2}^2 + \big[\lambda^{\text{opt}}_{\text{close}}(s)\big]^2
         -2\Re(a+bs^2)\lambda^{\text{opt}}_{\text{close}}(s)}
       {\abs{a+bs^2}^2 + \big[\lambda^{\text{opt}}_{\text{close}}(s)\big]^2
         +2\Re(a+bs^2)\lambda^{\text{opt}}_{\text{close}}(s)}.
\end{equation}
From this expression, it is clear that:
\begin{enumerate}
\item $\abs{\rho^\text{close}_{\lambda_\text{open}^{\text{OO2}}}(s)}^2\to1$
  if $\lambda^{\text{opt}}_{\text{close}}(s)\to\pm\infty$,
  $\lambda^{\text{opt}}_{\text{close}}(s)\to0$
  or $(a+bs^2)\to0$;
\item $\abs{\rho^\text{close}_{\lambda_\text{open}^{\text{OO2}}}(s)}^2<1$
  if $\lambda^{\text{opt}}_{\text{close}}(s)>0$ and $\Re(a+bs^2)>0$;
\item $\abs{\rho^\text{close}_{\lambda_\text{open}^{\text{OO2}}}(s)}^2>1$
  if $\lambda^{\text{opt}}_{\text{close}}(s)<0$ and $\Re(a+bs^2)>0$.
\end{enumerate}
As with EMDA, the two last results are opposed for $\Re(a+bs^2)<0$.
Furthermore, it is worth noticing that, since
$\lambda^{\text{opt}}_{\text{close}}$
is changing its sign more than twice
(at least for sufficiently large values of $k$),
$a$ and $b$ \emph{cannot} be optimized
to guarantee that
$\abs{\rho^\text{close}_{\lambda_\text{open}^{\text{OO2}}}(s)}^2<1$
$\forall{}s^2<k^2$.

\begin{figure}[ht]
  \centering
  \includegraphics{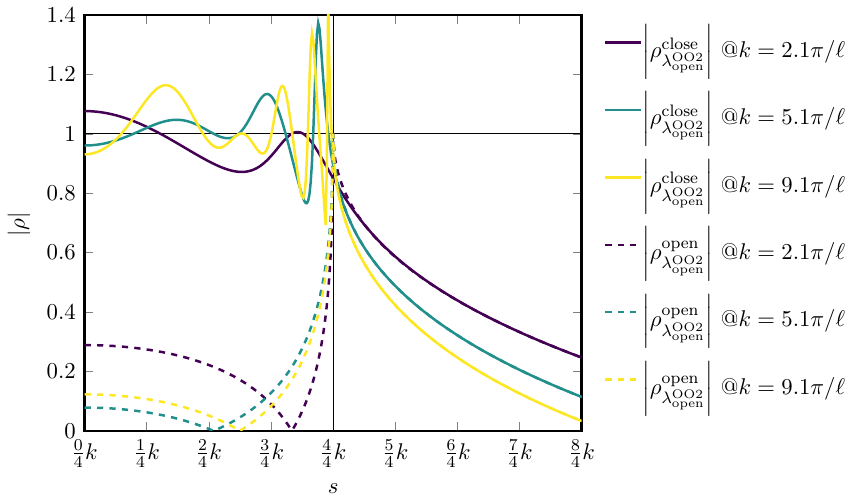}
  \caption{Graphs of
    $\abs{\rho^\text{close}_{\lambda_{\text{open}}^{\text{OO2}}}}$ and
    $\abs{\rho^\text{open}_{\lambda_{\text{open}}^{\text{OO2}}}}$
    for different values of $k$
    (the values the coefficients $a$ and $b$,
    as appearing in~\eqref{eq:open:lambda:oo2},
    are chosen according to~\cite{Gander2002}).}
  \label{fig:rho:oo2}
\end{figure}

Figure~\ref{fig:rho:oo2} compares
$\abs{\rho^\text{close}_{\lambda_\text{open}^{\text{OO2}}}(s)}$
with
$\abs{\rho^\text{open}_{\lambda_\text{open}^{\text{OO2}}}(s)}$
(\ie the convergence radius of the OO2 operator when used
in an OS scheme for solving the unbounded problem~\eqref{eq:helmholtz:open}).
Again, as expected, the distance between both radii
becomes negligible as $s$ grows (once $s^2>k^2$).
Nonetheless, when analyzing non-evanescent modes,
the performance of OO2 is poorer for cavity problems
than for unbounded ones since:
\begin{equation*}
  \max_s\abs{\rho^\text{open}_{\lambda_{\text{open}}^{\text{OO2}}}(s)}<1~%
  \text{(see~\cite{Gander2002})}
  \qquad
  \max_s\abs{\rho^\text{close}_{\lambda_{\text{open}}^{\text{OO2}}}(s)}>1
  \qquad
  \forall{}s^2<k^2.
\end{equation*}

\subsubsection{Square-root operator and its Pad\'e localization (PADE)}
In order to assess the performance of the PADE operator,
let us first determine the performance of the square-root
operator $\Trsm^{\text{opt}}_{\text{open}}$, as given in~\eqref{eq:open:S:opt},
when used in the OS scheme~\eqref{eq:ddm} solving
our model cavity problem~\eqref{eq:helmholtz}.
Indeed, as the Pad\'e localization process leads to an excellent approximation
of $\Trsm^{\text{opt}}_{\text{open}}$,
analyzing this limit case will shed light on the performance of the PADE operator,
at least for a sufficiently large number of Pad\'e terms $N_p$.

Given the expressions of $\lambda_{\text{open}}^{\text{opt}}$
in~\eqref{eq:open:lambda:opt}, we can deduce that the convergence radius
$\rho^\text{close}_{\lambda_{\text{open}}^{\text{opt}}}$
writes
\begin{align*}[left = {\rho^\text{close}_{\lambda_{\text{open}}^{\text{opt}}} = \empheqlbrace}]
  \frac
  {-\imag{}k\sqrt{1-\frac{s^2}{k^2}}-\sqrt{k^2-s^2}
    \cot\mathopen{}\left[
      \sqrt{k^2-s^2}\,\frac{\ell}{2}\mathclose{}
    \right]}
  {-\imag{}k\sqrt{1-\frac{s^2}{k^2}}+\sqrt{k^2-s^2}
    \cot\mathopen{}\left[
      \sqrt{k^2-s^2}\,\frac{\ell}{2}\mathclose{}
  \right]}
  & =
  \frac
  {-\imag-\cot\mathopen{}\left[
      \sqrt{k^2-s^2}\,\frac{\ell}{2}\mathclose{}
    \right]}
  {-\imag+\cot\mathopen{}\left[
      \sqrt{k^2-s^2}\,\frac{\ell}{2}\mathclose{}
  \right]}\hspace{-5pt}
  & \text{if}~s^2<k^2,\\
  \frac{0-2/\ell}{0+2/\ell}
  & = -1
  & \text{if}~s^2=k^2,\\
  \frac
  {k\sqrt{\frac{s^2}{k^2}-1}-\sqrt{s^2-k^2}
    \coth\mathopen{}\left[
      \sqrt{s^2-k^2}\,\frac{\ell}{2}\mathclose{}
    \right]}
  {k\sqrt{\frac{s^2}{k^2}-1}+\sqrt{s^2-k^2}
    \coth\mathopen{}\left[
      \sqrt{s^2-k^2}\,\frac{\ell}{2}\mathclose{}
  \right]}
  & =
  \frac
  {1-\coth\mathopen{}\left[
      \sqrt{s^2-k^2}\,\frac{\ell}{2}\mathclose{}
    \right]}
  {1+\coth\mathopen{}\left[
      \sqrt{s^2-k^2}\,\frac{\ell}{2}\mathclose{}
    \right]}
  & \text{if}~s^2>k^2,
\end{align*}
and therefore:
\begin{subequations}
  \label{eq:close:rho:opti}
  \begin{align}[left = {\abs{\rho^\text{close}_{\lambda_{\text{open}}^{\text{opt}}}} = \empheqlbrace}]
    & 1
    & \text{if}~s^2\leq{}k^2,\\
    & \exp\mathopen{}\left[-\ell\sqrt{s^2-k^2}\mathclose{}\right]
    & \text{if}~s^2>k^2,
  \end{align}
\end{subequations}
where the last line is obtained by directly exploiting the definition of
the hyperbolic cotangent~\cite{Oldham2009}.
Before studying the convergence radius~\eqref{eq:close:rho:opti},
it is worth mentioning that in the case $s^2<k^2$,
$\rho^\text{close}_{\lambda_{\text{open}}^{\text{opt}}}$
might be undefined because of the cotangent.
However, the limit for~$s$ approaching a pole of the cotangent is well defined
and is equal to $-1$.
This can be easily shown by a direct application of L'H\^opital's rule.
As expected, we have that
$\abs{\rho^\text{close}_{\lambda_{\text{open}}^{\text{opt}}}}$
is identically equal to $1$ for $s^2\leq{}k^2$ and decrease exponentially to $0$
(once $s^2>k^2$) as $s$ grows.

\begin{figure}[ht]
  \centering
  \includegraphics{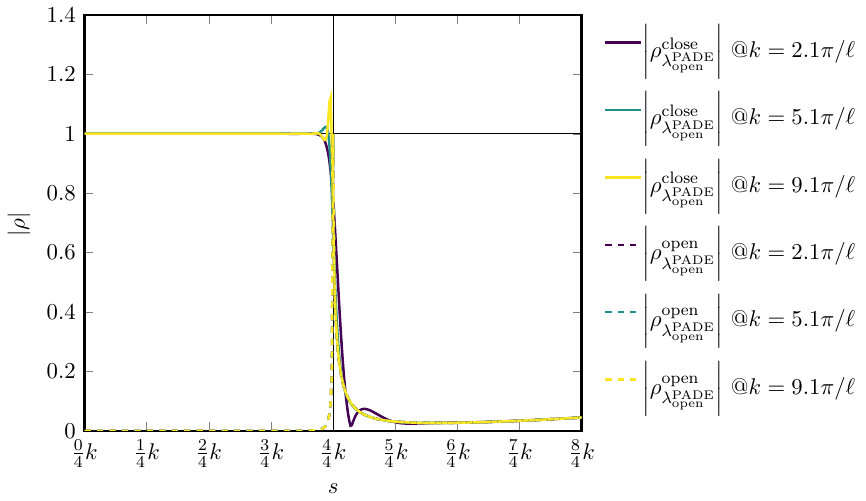}
  \caption{Graphs of
    $\abs{\rho^\text{close}_{\lambda_{\text{open}}^{\text{PADE}}}}$ and
    $\abs{\rho^\text{open}_{\lambda_{\text{open}}^{\text{PADE}}}}$
    for different values of $k$ an with $N_p=4$.}
  \label{fig:rho:pade}
\end{figure}

Let us now come back to the PADE operator.
As its symbol~\eqref{eq:open:lambda:pade}
is obtained from a Pad\'e approximation of $\lambda_{\text{open}}^{\text{opt}}$,
and as this approximation converges to $\lambda_{\text{open}}^{\text{opt}}$
as the number of Pad\'e terms $N_p\to\infty$~\cite{Antoine2006},
we can argue that the convergence radius determined in~\eqref{eq:close:rho:opti}
for the non-local square-root operator
approximates sharply the convergence radius
$\rho^\text{close}_{\lambda_{\text{open}}^{\text{PADE}}}$
of the PADE operator, at least for sufficiently large values of $N_p$,
as show in Figure~\ref{fig:rho:pade}.

\subsubsection{Summary}
Before concluding this section,
let us summarize our analysis.
When $s^2<k^2$, we showed that
$\abs{\rho^\text{close}_{\lambda_{\text{open}}^{\text{000, EMDA, OO2, PADE}}}}$
can be greater than $1$.
This behavior is significantly different from the unbounded case where
$\abs{\rho^\text{open}_{\lambda_{\text{open}}^{\text{000, EMDA, OO2, PADE}}}}$
is always smaller than $1$.
On the other hand, when $s^2>k^2$,
the difference between
$\abs{\rho^\text{close}_{\lambda_{\text{open}}^{\text{000, EMDA, OO2, PADE}}}}$
and
$\abs{\rho^\text{open}_{\lambda_{\text{open}}^{\text{000, EMDA, OO2, PADE}}}}$
vanishes exponentially.
Therefore, when comparing the OO2, EMDA, OO2 or PADE operators for solving
\begin{myenum*}
\item a cavity problem similar to~\eqref{eq:helmholtz} with
\item an unbounded problem without obstacles
  similar to~\eqref{eq:helmholtz:open},
\end{myenum*}
we can conclude that
the performance of the OS scheme will be deteriorated in the cavity case.

\section{Numerical experiments and curved geometries}
\label{sec:numexp}
In this section, we carry out some numerical experiments
in order to validate the previous theoretical analysis.
Furthermore, as we restricted ourselves to rectangular geometries,
we also offer in this section a numerical study of the performance
of the aforementioned OS schemes when handling curved configurations:
\ie circles (2D) and spheres (3D).

\subsection{Validation}
\begin{figure}[ht]
  \centering
  \includegraphics{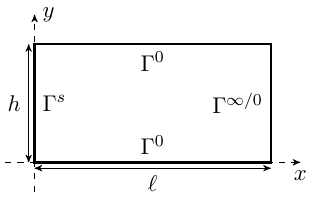}
  \caption{Rectangular geometry used for the numerical experiments.}
  \label{fig:numexp:geo:rect}
\end{figure}

Let us now illustrate the performance deterioration of
the OO0, EMDA, OO2 and PADE transmission conditions,
when used in an OS scheme for closed-domain Helmholtz problems.
To this end, two different cases will be presented:
\begin{myenum*}
\item a closed two-dimensional rectangular cavity with a length $\ell=9.5\lambda_w$,
  where $\lambda_w$ is the wavelength and
\item a section of an infinitely long two-dimensional rectangular waveguide
  with the same length $\ell$.
\end{myenum*}
In particular, the geometry displayed in Figure~\ref{fig:numexp:geo:rect}
is used and the following boundary conditions are imposed:
\begin{equation}
  \label{eq:numexp:bc}
  \begin{array}{r@{\,}l@{\quad}ll}
    p & = 0 & \text{on}~\Gamma^{\infty/0} & \text{in the cavity case},\\
    \vec{n}\cdot\Grad{p} & = \imag{}k
            & \text{on}~\Gamma^{\infty/0} & \text{in the waveguide case},\\
    p & = 0 & \text{on}~\Gamma^0 & \text{in both cases},\\
    p & = \displaystyle\sum_{m=1}^N\sin\left(
        m\frac{\pi}{h}y\right) & \text{on}~\Gamma^s & \text{in both cases},\\
  \end{array}
\end{equation}
where $k=\frac{2\pi}{\lambda_w}$, $h=\ell/2$ and $N=9$.
Let us stress that, as the length of the cavity is not an integer multiple
of the wavelength, the closed Helmholtz problem is well defined.
Furthermore, let us also note that $9$ modes can propagate in the waveguide
for the selected height and wavenumber.
All these modes are superimposed
when exciting both the cavity and the waveguide problems.

Concerning the numerical setup,
the Helmholtz problem~\eqref{eq:helmholtz} is discretized
with a finite element method of order $5$
and the geometry in Figure~\ref{fig:numexp:geo:rect} is discretized
with $8$ triangular mesh elements per wavelength.
An optimized Schwarz scheme is then used to solve~\eqref{eq:helmholtz:omega}
with the boundary conditions given in~\eqref{eq:numexp:bc}, with $g=0$
and with two subdomains of equal size,
as shown in Figure~\ref{fig:domain}.
Furthermore, let us mention that in the following numerical experiment,
the non-overlapping fixed-point Schwarz algorithm is recast
into the linear system~\cite{Dolean2015a}:
\begin{equation}
  \label{eq:numexp:lin}
  (\Ident-\OpA)\vec{d} = \vec{b},
\end{equation}
where one application of the operator $\OpA$ amounts to one iteration of the
fixed-point procedure with \emph{homogeneous} Dirichlet boundary conditions,
where $\Ident$ is the identity operator,
where the vector $\vec{d}$ concatenates all $\vec{n}\cdot\Grad{p}+\Trsm(p)$
at the interface between the subdomains
and where the right hand side vector $\vec{b}$ results from the
\emph{non-homogeneous} Dirichlet boundary conditions.
This linear system is then solved with a matrix-free GMRES
\emph{without} restart.
Regarding the free parameters of the transmission conditions, we choose:
$4$ Pad\'e terms and no damping for the PADE transmission condition~\cite{Boubendir2012};
the parameters $\alpha, \beta$ for the OO2 operator
as proposed in~\cite{Gander2002}
and a damping of $25\%$ for the EMDA operator
as suggested in~\cite{Marsic2016b}.

\subsubsection{Convergence of GMRES}
\begin{figure}[ht]
  \centering
  \includegraphics{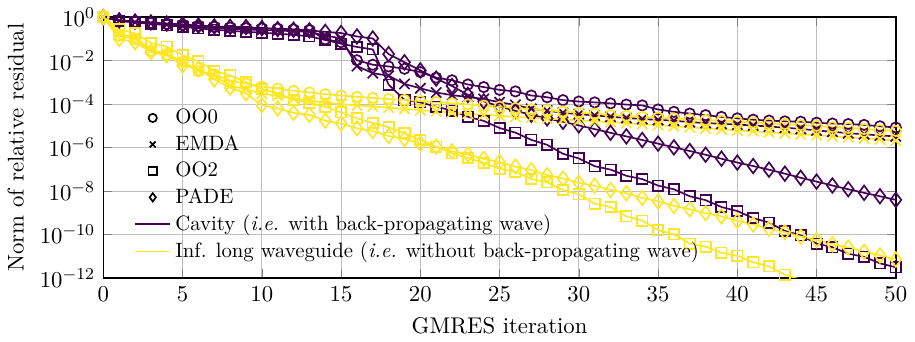}
  \caption{Convergence profile of GMRES
    for an infinitely long waveguide and for a cavity
    (two-dimensional rectangular case).}
  \label{fig:numexp:gmres}
\end{figure}

As a first numerical experiment, let us analyze the convergence rate
of GMRES for solving both cavity and waveguide problems,
as displayed in Figure~\ref{fig:numexp:gmres}.
From these data, the performance loss in case of back-propagating waves
(\ie the cavity scenario) is obvious:
a difference of 2 orders of magnitudes in the relative GMRES residual
(between both scenarios) for OO2 and PADE.

\subsubsection{Spectrum of $\Ident - \OpA$}
\begin{figure}[ht]
  \centering
  \includegraphics[width=7cm]{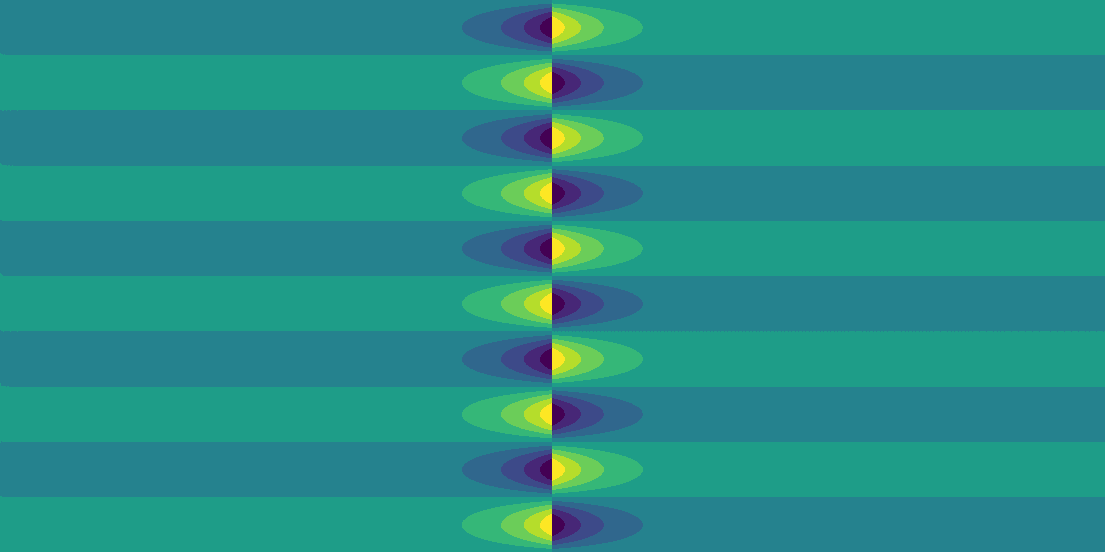}\hspace{1cm}
  \includegraphics[width=7cm]{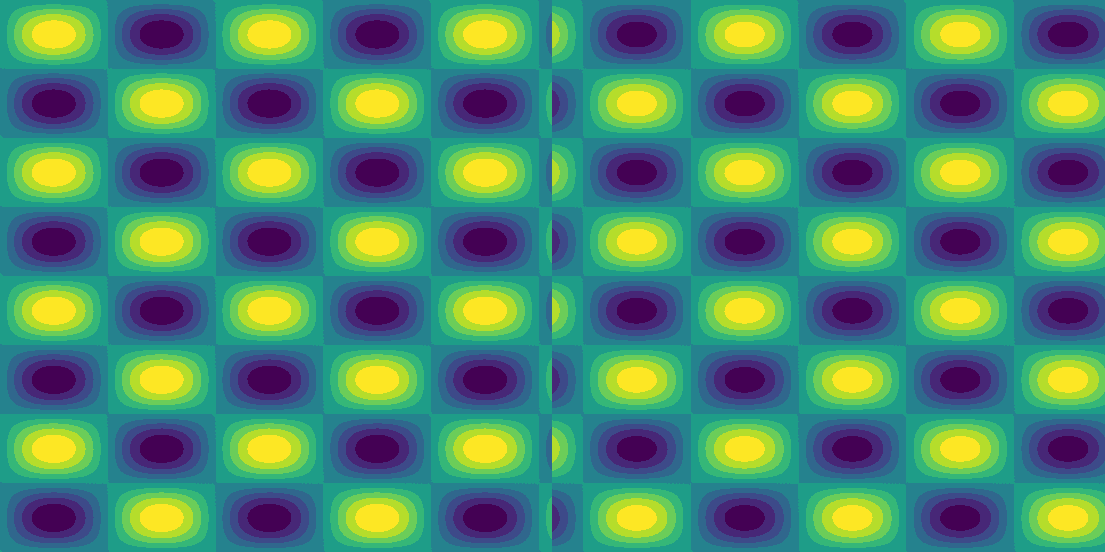}
  \caption{Evanescent (left) and non-evanescent (right) eigenmodes
    of $(\Ident-\OpA)$
    (two-dimensional rectangular case).}
  \label{fig:numexp:eigV}
\end{figure}

\begin{figure}[p]
  \centering
  \begin{tabular}{cc}
    \includegraphics{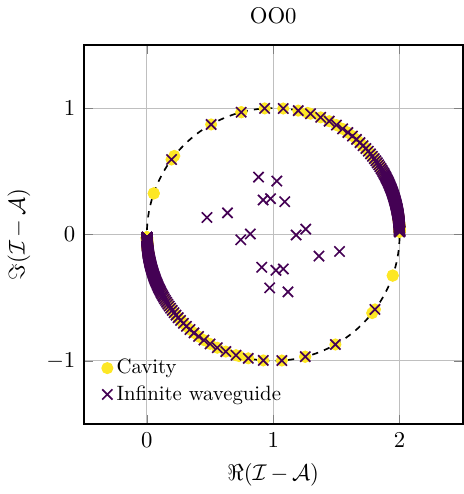} &
    \includegraphics{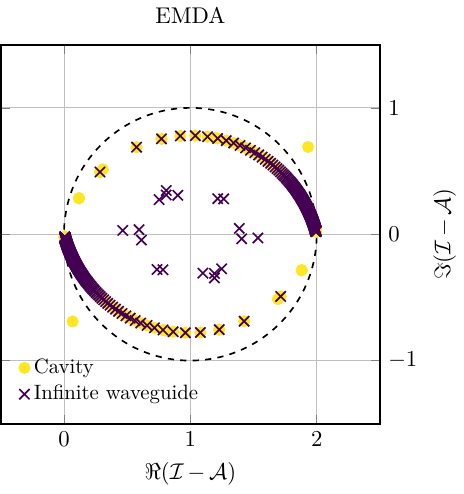}\\
    \includegraphics{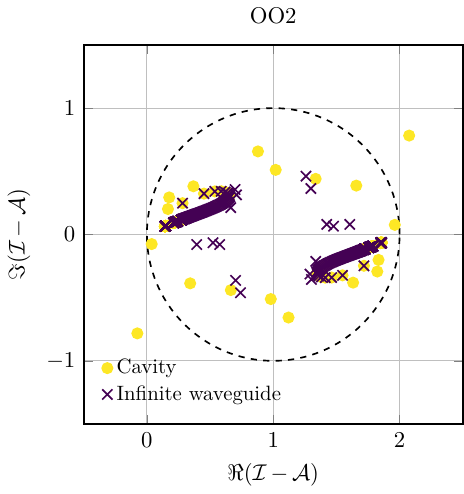} &
    \includegraphics{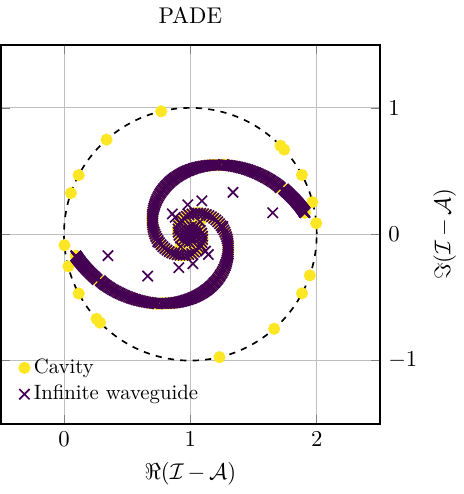}
  \end{tabular}
  \caption{Spectrum of $\Ident - \OpA$ with different transmission conditions
    (two-dimensional rectangular case).}
  \label{fig:numexp:spect:rect}
\end{figure}

In a second numerical experiment,
the spectrum of the system matrix $\Ident - \OpA$
is studied for both rectangular cavity and waveguide cases.
As shown in Figure~\ref{fig:numexp:spect:rect},
and as predicted by the theory,
we can observe that for \emph{cavity} problems:
\begin{enumerate}
\item all the eigenvalues lie \emph{on} the unit circle
  for OO0;
\item some eigenvalues are located \emph{outside} the unit circle
  for OO2 and EMDA and
\item all $9$ non-evanescent modes lie \emph{on} the unit circle
  for PADE.
\end{enumerate}
For illustration purposes, a non-evanescent mode and an evanescent one
are displayed in Figure~\ref{fig:numexp:eigV}.
These modes are nothing else but eigenvectors of $(\Ident-\OpA)$.

\clearpage
\subsection{Curved geometries}
In the second part of this section, we study
the behavior of the OO0, EMDA, OO2 and PADE transmission conditions
when used in curved geometries (either open or closed) as shown
in Figures~\ref{fig:numexp:geo:cir} and~\ref{fig:numexp:geo:sph} by numerical experiments.
Let us note that in order to simplify the discussion,
this subsection will consider \emph{source-free} problems
and thus focus on the spectrum of $\Ident - \OpA$ only.

\subsubsection{Two-dimensional case: a circular geometry}
\begin{figure}[ht]
  \centering
  \includegraphics{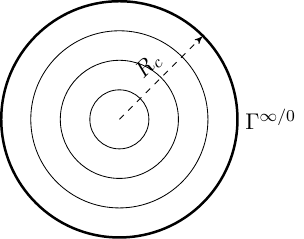}
  \caption{Two-dimensional circular geometry
    used for the numerical experiments.}
  \label{fig:numexp:geo:cir}
\end{figure}

Let us start by considering a two-dimensional circular geometry
with a radius $R_c$ such that $2R_c = 9.5\lambda_w$.
In contrast to the previous test case,
a circular (concentric) domain decomposition is chosen,
and
in order to consider a configuration with more than two subdomains
we adopt here a \emph{four} subdomains setting.
Regarding the discretization,
we use again the boundary conditions~\eqref{eq:numexp:bc}\footnote{Let us note
  that the boundary conditions on $\Gamma^0$ and $\Gamma^s$ are simply ignored
  in this case.
  Furthermore, the condition $\vec{n}\cdot\Grad{p} = \imag{}k$
  in~\eqref{eq:numexp:bc}
  is no longer an \emph{exact} absorbing boundary condition
  in this curved setting.}
and discretize the Helmholtz problem with a finite element method of order $5$
on a geometry meshed with $8$ \emph{curved} (second-order) triangular elements
per wavelength.
Concerning the free parameters of the transmission operators, we choose:
4 Pad\'e terms without damping for the PADE condition;
the parameters $\alpha$, $\beta$ of the OO2 condition
as proposed in~\cite{Boubendir2012}
and a damping of $50\%$ for the EMDA condition
as suggested in~\cite{Boubendir2007}.
For illustration purposes, a schematic representation of the numerical setting
is shown in Figure~\ref{fig:numexp:geo:cir}.

\begin{figure}[p]
  \centering
  \begin{tabular}{cc}
    \includegraphics{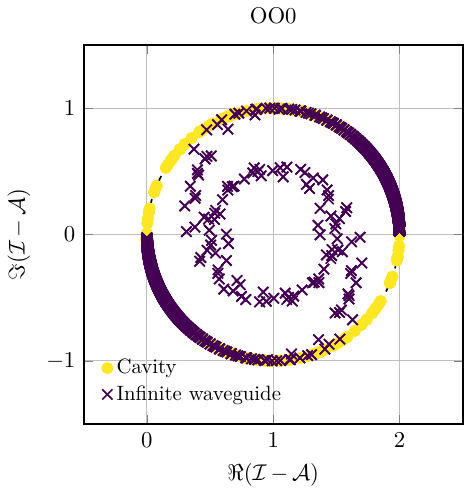} &
    \includegraphics{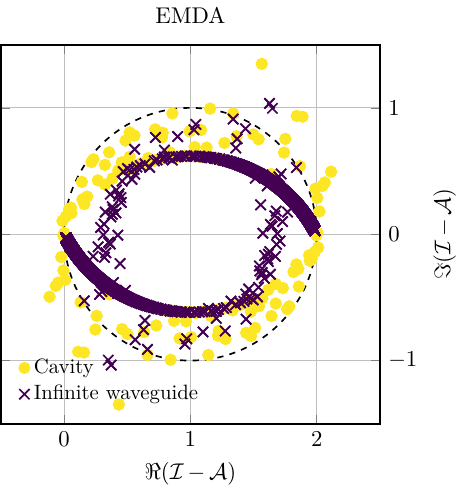}\\
    \includegraphics{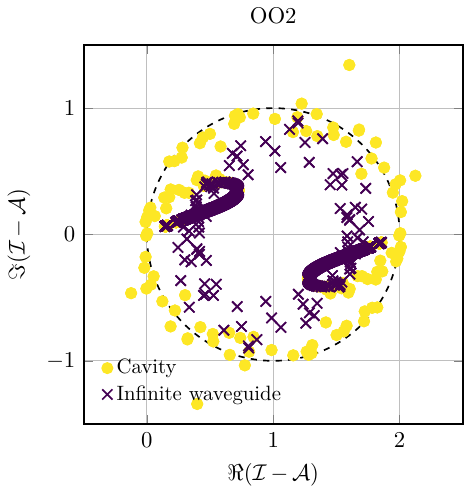} &
    \includegraphics{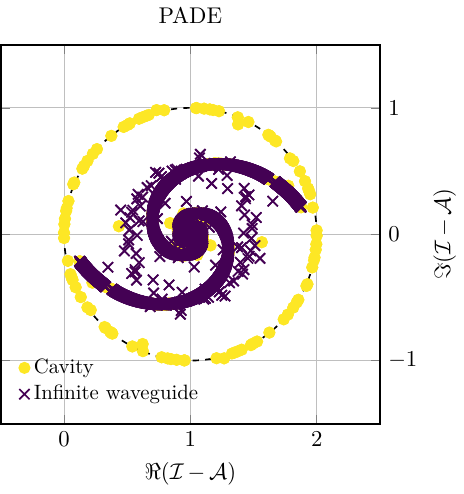}
  \end{tabular}
  \caption{Spectrum of $\Ident - \OpA$ with different transmission conditions
    (two-dimensional circular case with four subdomains).}
  \label{fig:numexp:spect:cir}
\end{figure}

\begin{figure}[ht]
  \centering
  \includegraphics{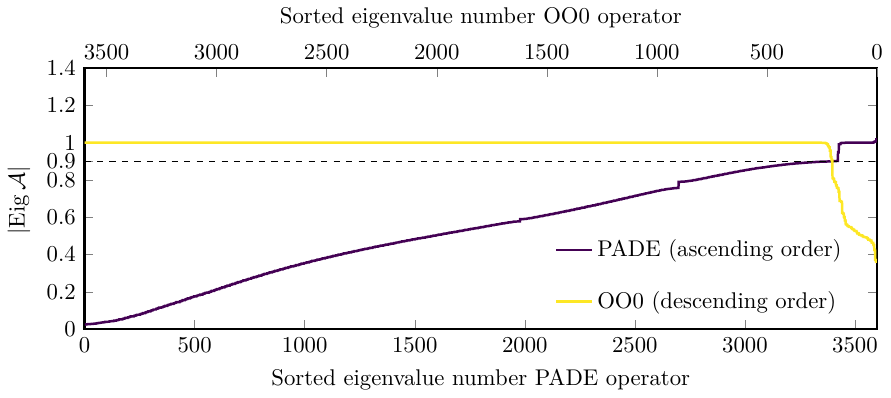}
  \caption{Modulus of the eigenvalues of $\OpA$
    with PADE (closed problem) and OO0 (open problem);
    the eigenvalues are sorted in:
    \textit{i)}~ascending order with PADE and
    \textit{ii)}~descending order with OO0
    (two-dimensional circular case with four subdomains).}
  \label{fig:numexp:rho:cir}
\end{figure}

The spectra of $\Ident - \OpA$ for open and closed problems
and for the different transmission conditions are depicted in Figure~\ref{fig:numexp:spect:cir}.
As it can be directly seen from those data,
the conclusions drawn in section~\ref{sec:open:S} are once again recovered.
Indeed, in the case of a cavity problem:
\begin{myenum*}
\item all eigenvalues lie on the unit circle with the OO0 condition;
\item some eigenvalues are outside the unit circle
  with the EMDA and OO2 operators and
\item some eigenvalues lie on the unit circle with the PADE condition.
\end{myenum*}
Let us also note that the number of eigenvalues lying
on the unit circle for the PADE operator in the closed case
should match the number of eigenvalues
\emph{inside} the unit circle for the OO0 condition in the \emph{open} case.
This behavior is recovered,
up to a few modes near the unit circle (\ie $\abs{\Eig{\OpA}}\simeq0.9$),
as shown in Figure~\ref{fig:numexp:rho:cir}.

\clearpage

\subsubsection{Three-dimensional case: a spherical geometry}
\begin{figure}[ht]
  \centering
  \includegraphics{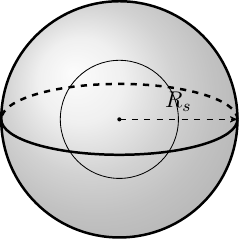}
  \caption{Three-dimensional spherical geometry
    used for the numerical experiments.}
  \label{fig:numexp:geo:sph}
\end{figure}

For this last numerical experiment, a three-dimensional case is studied.
In particular, we extend the previous two-dimensional circular scenario to a
three-dimensional spherical problem, as shown in Figure~\ref{fig:numexp:geo:sph}.
However, in order to keep the computational effort reasonable,
this last numerical experiment is restricted to:
\begin{myenum*}
\item a sphere radius $R_s$ such that $2R_s = 3.5\lambda_w$;
\item two subdomains and
\item a second-order finite element discretization.
\end{myenum*}
All other numerical aspects,
such as the free parameters of the transmission operators,
are kept identical to the previous study.

\begin{figure}[p]
  \centering
  \begin{tabular}{cc}
    \includegraphics{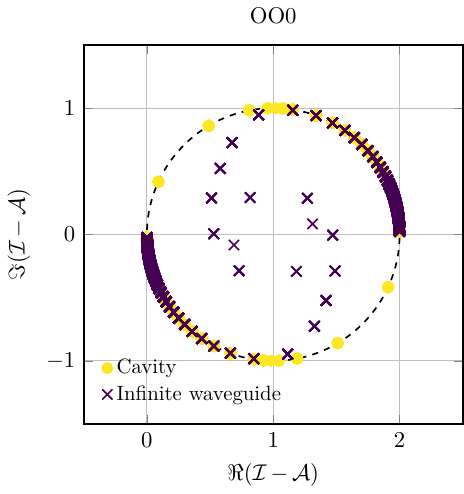} &
    \includegraphics{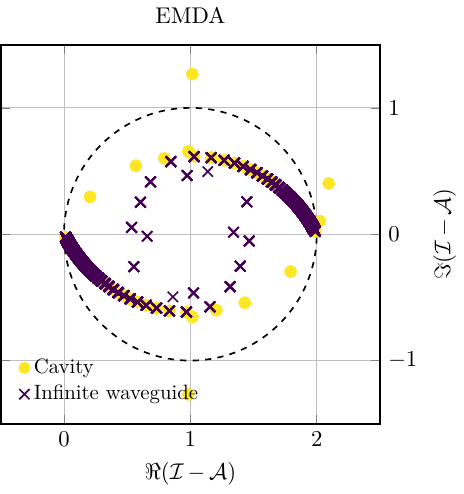}\\
    \includegraphics{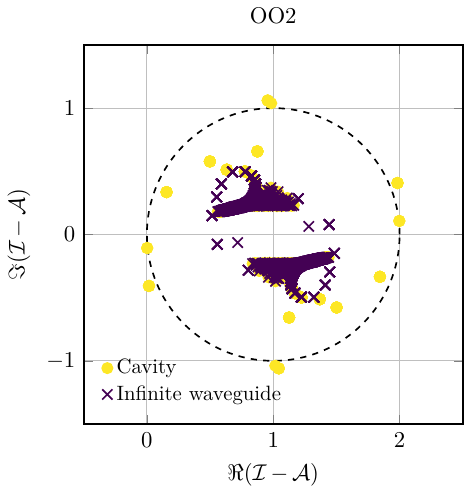} &
    \includegraphics{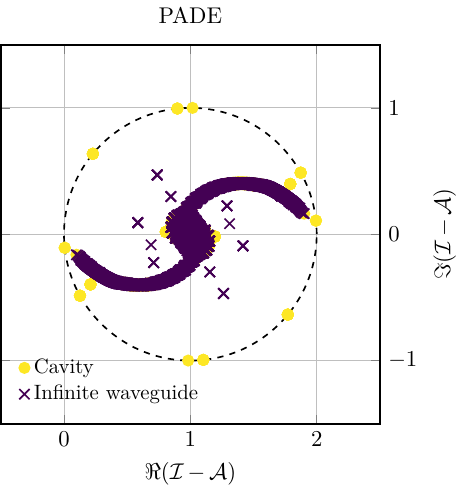}
  \end{tabular}
  \caption{Spectrum of $\Ident - \OpA$ with different transmission conditions
    (three-dimensional spherical case).}
  \label{fig:numexp:spect:sph}
\end{figure}

\begin{figure}[ht]
  \centering
  \includegraphics{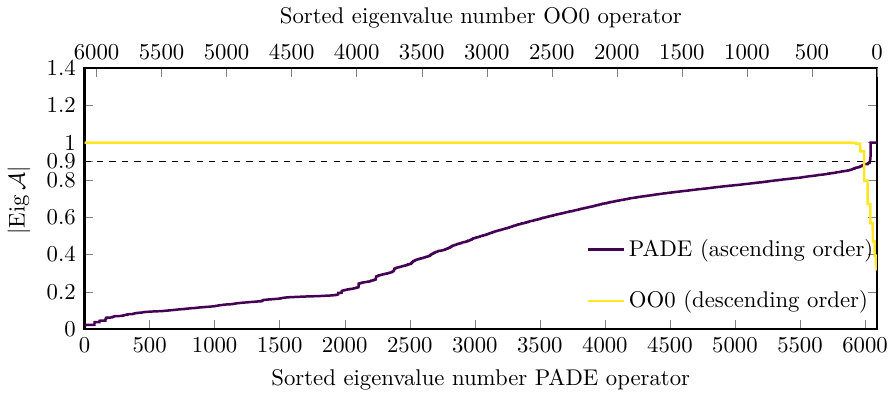}
  \caption{Modulus of the eigenvalues of $\OpA$
    with PADE (closed problem) and OO0 (open problem);
    the eigenvalues are sorted in:
    \textit{i)}~ascending order with PADE and
    \textit{ii)}~descending order with OO0
    (three-dimensional spherical case with two subdomains).}
  \label{fig:numexp:rho:sph}
\end{figure}

Once again, we compute the spectra of $\Ident - \OpA$
for open and closed problems
and for the different transmission conditions.
From the data depicted in Figure~\ref{fig:numexp:spect:sph},
it is obvious that the conclusions we have drawn
in the previous two-dimensional numerical experiments remain once again valid.
In particular:
\begin{myenum*}
\item all eigenvalues lie on the unit circle with the OO0 condition;
\item some eigenvalues are outside the unit circle
  with the EMDA and OO2 operators and
\item some eigenvalues lie on the unit circle with the PADE condition.
\end{myenum*}
Regarding this last point, we recover that the number of eigenvalues
lying on the unit circle for the PADE operator in the closed case matches
(up to a few eigenvalues near $\abs{\Eig{\OpA}}\simeq0.9$)
the number of eigenvalues \emph{inside} the unit circle for the OO0 condition
in the \emph{open} case, as shown in Figure~\ref{fig:numexp:rho:sph}.

\section{Conclusion}
\label{sec:conclusion}
In this paper, we derived the optimal transmission operator
$\Trsm^{\text{opt}}_{\text{close}}$ of an optimized Schwarz scheme
solving a simple rectangular Helmholtz cavity problem.
We furthermore demonstrated that for \emph{evanescent} modes
this optimal operator is excellently approximated by
the optimal transmission operator of open problems without obstacles
$\Trsm^{\text{opt}}_{\text{open}}$.
On the other hand, we also showed that $\Trsm^{\text{opt}}_{\text{open}}$
cannot be used to approximate $\Trsm^{\text{opt}}_{\text{close}}$
for \emph{non-evanescent} modes.
For this reason the classical OO0, EMDA, OO2 and PADE transmission operators
exhibit a performance drop when applied to cavity problems
(compared to an equivalent unbounded configuration).
In particular, we determined that the convergence radius
$\rho^\text{close}(s)$ of the OS scheme exhibits:
\begin{enumerate}
\item a modulus equal to $1$ for all $s\in\mathbb{R}$ for the OO0 operator;
\item a modulus greater than $1$ for some $s^2<k^2$
  for the EMDA and OO2 operators and
\item a modulus equal to $1$ for all $s^2\leq{}k^2$ for the PADE operator,
  at least when using a sufficiently large number of Pad\'e terms.
\end{enumerate}
Furthermore, by means of numerical experiments,
we showed that those conclusions,
which were derived for a rectangular domain,
remain valid when considering curved and three-dimensional configurations.

\section*{Acknowledgments}
This research project has been funded by the Deutsche Forschungsgemeinschaft
(DFG, German Research Foundation) -- Project number 445906998.
The authors would like to express their gratitude to
Ms.~Heike Koch, Mr.~Achim Wagner, Mr.~Dragos Munteanu and Dr.~Wolfgang M\"uller
for the administrative and technical support.

\bibliographystyle{ieeetr}
\bibliography{ddm.bib}
\end{document}